\documentclass[11pt]{amsart}
\usepackage[T1]{fontenc}
\usepackage{times}
\usepackage{amsmath,amssymb,amsthm,mathtools,mathrsfs}
\usepackage{microtype,enumitem,esint}
\usepackage{graphicx,xcolor,tikz,placeins,float}
\usetikzlibrary{arrows.meta,calc,positioning,fit,patterns,backgrounds}
\usepackage[colorlinks=true,linkcolor=red,citecolor=red,urlcolor=blue,
  hypertexnames=false,bookmarksnumbered=true]{hyperref}
\usepackage{bookmark}
\usepackage{accents} 

\allowdisplaybreaks[2]
\setlist[enumerate]{label=\textup{(\roman*)},leftmargin=*,itemsep=2pt,topsep=4pt}
\theoremstyle{plain}
\newtheorem{theorem}{Theorem}[section]
\newtheorem{lemma}[theorem]{Lemma}

\newtheorem{corollary}[theorem]{Corollary}
\theoremstyle{definition}
\newtheorem{definition}[theorem]{Definition}

\numberwithin{equation}{section}
\newcommand{\R}{\mathbb R}
\newcommand{\N}{\mathbb N}

\newcommand{\dist}{\operatorname{dist}}

\newcommand{\divg}{\operatorname{div}}

\newcommand{\dd}{\,d}
\newcommand{\avint}{\fint}
\newcommand{\Hn}{\R^{n+1}_+}
\DeclareMathOperator{\esssup}{ess\,sup}
\definecolor{mathBlue}{RGB}{34,79,119}
\definecolor{mathOrange}{RGB}{190,98,28}
\definecolor{mathGreen}{RGB}{35,125,100}
\definecolor{lightFill}{RGB}{239,244,248}

\newcommand{\DRAffineFigure}{%
\begin{figure}[htbp]
\centering
\begin{tikzpicture}[font=\small,>=Latex,line cap=round,line join=round]
 \begin{scope}[shift={(-4.0,-.15)},
   x={(.85cm,-.14cm)},y={(.48cm,.28cm)},z={(0cm,.82cm)}]
  \path[fill=lightFill] (0,0,0)--(3,0,0)--(3,3,0)--(0,3,0)--cycle;
  \draw[mathBlue!45] (0,0,0)--(3,0,0)--(3,3,0)--(0,3,0)--cycle;
  \draw[->,mathBlue,line width=.9pt] (0,0,0)--(3.5,0,0)
      node[right] {$e_1$};
  \draw[->,mathBlue,line width=.9pt] (0,0,0)--(0,3.5,0)
      node[above right] {$e_2$};
  \draw[->,mathGreen,line width=1pt] (0,0,0)--(0,0,2.8)
      node[above] {$e_3$};
  \draw[<->,mathOrange,line width=1pt]
      (2.2,0,0) to[bend right=24] (0,2.2,0);
  \node[mathGreen,anchor=east] at (-.1,0,1.8) {$Se_3=0$};
  \node[anchor=north] at (2.0,1.5,-.7) {$t=0$};
 \end{scope}
 \node[align=left,anchor=west,text width=7.4cm] at (1.35,1.55)
  {$S=e_1\otimes e_2+e_2\otimes e_1$,\quad $A=I+b(x_2,t)S$,\\[6pt]
   $u_0(x)=x_1+v(x_2,t)$,\quad $Ae_3=e_3$.};
 \node[align=left,anchor=west,text width=7.4cm] at (1.35,.10)
  {$\displaystyle A\nabla u_0=(1+bv_{x_2},\,b+v_{x_2},\,v_t)$,\\[6pt]
   $\displaystyle \operatorname{div}(A\nabla u_0)
     =b_{x_2}+v_{x_2x_2}+v_{tt}$.};
\end{tikzpicture}
\caption{The equation for $u_0=x_1+v(x_2,t)$ in $\mathbb R^3$. The matrix changes only
the two tangential coordinates. The height function remains a solution,
while the correction to $x_1$ satisfies a scalar equation in $(x_2,t)$.
The orange arrows indicate $Se_1=e_2$ and $Se_2=e_1$.}
\label{fig:dr-affine}
\end{figure}%
}

\newcommand{\DRObservationFigure}{%
\begin{figure}[htbp]
\centering
\begin{tikzpicture}[font=\small,>=Latex,line cap=round,line join=round]
 \begin{scope}[shift={(-4.2,0)}]
  \path[fill=lightFill] (0,0) rectangle (5.0,3.25);
  \foreach \lo/\hi/\lab in {2.65/2.98/N-2,1.95/2.28/N-1,1.25/1.58/N}{
   \path[fill=mathBlue!20,draw=mathBlue!65]
       (0,\lo) rectangle (5.0,\hi);
   \node[anchor=west] at (5.12,{.5*(\lo+\hi)}) {$i=\lab$};
  }
  \path[fill=mathGreen!16,draw=mathGreen!65] (0,.45) rectangle (5.0,.60);
  \path[fill=mathGreen!16,draw=mathGreen!65] (0,.16) rectangle (5.0,.31);
  \node[anchor=west,mathGreen] at (5.12,.52) {$i=N+1$};
  \node[anchor=west,mathGreen] at (5.12,.23) {$i=N+2$};
  \draw[mathBlue,line width=.8pt] (0,0)--(5.0,0);
  \draw[->] (-.3,0)--(-.3,3.5) node[above] {$t$};
  \draw[->] (0,-.18)--(1.0,-.18) node[right] {$x_2$};
  \draw[mathOrange,dashed] (0,.95)--(4.9,.95);
  \draw[mathOrange,line width=.8pt] (2.50,.95) circle (.18);
  \fill[mathOrange] (2.50,.95) circle (1.7pt);
  \node[anchor=east,mathOrange] at (-.42,.95) {$t_N$};
  \draw[mathOrange,->] (3.30,-.55)--(2.60,.79);
  \node[anchor=west,mathOrange] at (1.60,-.72) {$B(x_N,t_N/2)$};
 \end{scope}
 \node[anchor=north west,align=left] at (3.1,3.25)
  {$\lambda_i=4^{m_i}$,\qquad $t_N=b_N/4$.};
 \node[anchor=north west,align=left] at (3.1,2.45)
  {$\displaystyle S_N(y)=\sum_{i\leq N}a_ic_{L_i}\sin(\lambda_iy)$.};
 \node[anchor=north west,align=left] at (3.1,1.20)
  {$\displaystyle
   \left|\partial_tu_0(x)+S_N(\xi_2)\right|\leq C\varepsilon$,\\[8pt]
   $x\in B(x_N,t_N/2)$,\quad $x_N=(\xi',t_N)$.};
\end{tikzpicture}
\caption{The observation height and the coefficient layers. The
terms with $i\leq N$ give the partial sum $S_N$ in the normal
derivative. The remaining terms and the error in replacing the first
$N$ profiles by their boundary values are uniformly bounded. Heights
and layer thicknesses are schematic.}
\label{fig:dr-observation}
\end{figure}%
}

\newcommand{\Di}{\operatorname{D}}
\newcommand{\dv}{\operatorname{div}}
\newcommand{\loc}{\operatorname{loc}}   
\newcommand{\Reg}{\operatorname{R}}

\newcommand{\m}[1]{\mathcal{#1}}
  
\newcommand{\bb}[1]{\mathbb{#1}}
\newcommand{\wt}{\widetilde}
\newcommand{\ep}{\varepsilon}

\def\Xint#1{\mathchoice
	{\XXint\displaystyle\textstyle{#1}}%
	{\XXint\textstyle\scriptstyle{#1}}%
	{\XXint\scriptstyle\scriptscriptstyle{#1}}%
	{\XXint\scriptscriptstyle\scriptscriptstyle{#1}}%
	\!\int}
\def\XXint#1#2#3{{\setbox0=\hbox{$#1{#2#3}{\int}$}
		\vcenter{\hbox{$#2#3$}}\kern-.5\wd0}}

\def\dashint{\Xint-}

\def\YYint#1#2#3{{\setbox0=\hbox{$#1{#2#3}{\iint}$}
		\vcenter{\hbox{$#2#3$}}\kern-.51\wd0}}
 
\newcommand*{\dt}[1]{%
	\accentset{\mbox{\Large\bfseries .}}{#1}}

\hypersetup{pdftitle={Dirichlet solvability, regularity and connectivity for elliptic operators}}
\begin{document}
	
	\author[A. Merlo]{Andrea Merlo}
	\address{Andrea Merlo\\ Departamento de Matem\'aticas, Universidad del Pa\'is Vasco (UPV/EHU), Barrio Sarriena s/n 48940 Leioa, Spain and\\
		Ikerbasque, Basque Foundation for Science, Bilbao, Spain.}
	\email{andrea.merlo@ehu.eus}
	
	\author[M. Mourgoglou]{Mihalis Mourgoglou}
	\email[]{michail.mourgoglou@ehu.eus}
	\address{Departamento de Matemáticas, Universidad del País Vasco, UPV/EHU, Barrio Sarriena S/N 48940 Leioa, Spain and, IKERBASQUE, Basque Foundation for Science, Bilbao, Spain}
	
	\author[B. Poggi]{Bruno Poggi}
	\email[]{brunopoggi@pitt.edu}
	\address{Department of Mathematics, University of Pittsburgh, Pittsburgh, PA, 15213} 
	
	\title[Dirichlet does not imply Regularity]{Dirichlet does not imply Regularity for general elliptic operators in smooth domains}
	\date{\today}
	
	\thanks{A.M. was supported by the European Research Council through the ERC Starting Grant DARTS (grant agreement No. 101299428). A.M. was also supported by IKERBASQUE, and partially supported by the project PID2024-157724NB-I00, funded by MICIU/AEI/10.13039/501100011033/FEDER, EU, by the Europa Excelencia grant EUR2025-165042, funded by the Ministerio de Ciencia, Innovación y Universidades and Agencia Estatal de Investigación, and by  IT1875-26 (Basque Government). M.M. was supported  by IKERBASQUE and partially supported by the grant PID2024-157724NB-I00: of the Ministerio de Econom\'ia y Competitividad (Spain), and by  IT-1615-22 (Basque Government). B.P. was partially supported by an award from the Simons Foundation in the Travel Support for Mathematicians program, [SFI-MPS-TSM-00025559, BP]}
	
\begin{abstract} In the upper half-space $\mathbb R^{n+1}_+$, $n\geq2$, we show that there exists a real, diagonal, strongly elliptic matrix $A$ of bounded measurable coefficients, smooth in the interior of $\mathbb R^{n+1}_+$ and continuous on its closure, such that for the elliptic operator $L=-\operatorname{div}A\nabla$, the Dirichlet problem with data in $L^p$ is solvable, for every $p\in(1,\infty)$, but the Dirichlet-Regularity problem is not solvable for any $p\in(1,\infty)$. The counterexample has the following properties, in particular:
	\begin{itemize}
		\item $A$ equals the identity matrix at the boundary and their $L^\infty$ distance in $\mathbb R^{n+1}_+$ can be made arbitrarily small.
		\item $A$ satisfies that $t|\partial_t A|^2$ is a Carleson measure.
		\item $A$ fails to satisfy the Dahlberg-Kenig-Pipher  condition only slightly, in the sense that for every $q>2$, $t^q|\nabla A|^q\frac{dxdt}{t}$ is a vanishing Carleson measure, but $t|\nabla A|^2$ is not a Carleson measure.		
	\end{itemize} 
	The counterexample resolves a long-standing question of C. Kenig.
	
\end{abstract}
\maketitle
\begingroup\footnotesize\tableofcontents\endgroup 

\section{Introduction}
\label{sec:introduction}

The $L^p$-solvability of homogeneous boundary value problems for second-order linear  elliptic partial differential equations   with non-smooth coefficients  has been an area of intense and fruitful research in the last few decades. In this manuscript, we consider the validity of a long-standing question in the field \cite{kbook, kp2}; namely, whether on a smooth domain (henceforth we consider the upper-half space $\Omega=\bb R^{n+1}_+$, $n\geq2$), and given a generic, real, symmetric, strongly elliptic operator 
\[
L=-\dv A\nabla
\]
with bounded measurable coefficients, the solvability of the Dirichlet problem in $L^p$ on $\Omega$ for some $p\in(1,\infty)$ implies the solvability of the Dirichlet-Regularity problem in $L^{q}$ for some $q\in(1,\infty)$. Here, we construct a counterexample that shows that the answer is \emph{no}; and moreover, the counterexample has several interesting properties, which we shall detail further below.

First, let us define the salient concepts at hand. Throughout this article we take $\Omega=\bb R^{n+1}_+=\bb R^n\times(0,\infty)$, $n\geq2$. and identify $\partial\Omega$ with $\R^n$. Using $t$ as the transversal variable, we shall write $X=(x,t)$, $\delta(X):=\dist(X,\partial\Omega)=t$, $e_t=e_{n+1}$, and
$\sigma=\m H^n|_{\partial\Omega}$. We consider a  matrix function $A=(A_{ij})$ defined on $\overline\Omega$ which is always a real,   symmetric $(n+1)\times(n+1)$ matrix of    bounded measurable coefficients in $\Omega$ verifying the strong ellipticity conditions 
\begin{equation}\label{eq.elliptic}
	\lambda|\xi|^2\leq\sum_{i,j=1}^{n+1}A_{ij}(x)\xi_i\xi_j,\qquad \Vert A\Vert_{L^{\infty}(\Omega)}\leq\frac1{\lambda},\quad x\in\Omega,\quad\xi\in\bb R^{n+1},
\end{equation}  
for some $\lambda\in(0,1)$. We denote $L_A=-\dv A\nabla$, and omit the subscript if the matrix is understood from context.

For  $\xi\in\partial\Omega$, we define the \emph{non-tangential cone with vertex $\xi$} by
\begin{equation}\label{eq.cone}\nonumber
	\gamma(\xi):=\{X=(x,t)\in \Omega:|X-\xi|<2t\},
\end{equation}
Define the \emph{non-tangential maximal function} of  $u\in L^{\infty}_{\loc}(\Omega)$ by 
\[
\m N(u)(\xi):=\sup_{X\in\gamma(\xi)}|u(X)|,\quad\text{for }\xi\in\partial\Omega.
\]
Following \cite{kp}, we introduce the \emph{modified non-tangential maximal function} $\widetilde{\m N}$: for any $u\in L^2_{\loc}(\Omega)$, we write
\begin{equation}\label{eq.mnt}\nonumber
	\widetilde{\m N}(u)(\xi):=\sup_{X\in\gamma(\xi)}\Big(\dashint_{B(X, ~\frac t2)}|u(Y)|^2\,dY\Big)^{1/2},\qquad\xi\in\partial\Omega.
\end{equation}
 
\begin{definition}[The $L^p$ Dirichlet and regularity problems]\label{def.dirichlet} Let $p\in(1,\infty)$. We say that the  \emph{Dirichlet problem   for the operator $L$ with $L^{p}$ data   is solvable} (write $(\Di_{p}^L)$ is solvable), if there exists $C\geq1$ so that for each $f\in C_c(\partial\Omega)$, the solution $u$ to the continuous Dirichlet problem
	\begin{equation}\label{eq.dirichlet}
		\left\{\begin{array}{ll}-\dv A\nabla u=0,\quad&\text{in }\Omega,\\u=f,\quad&\text{on }\partial\Omega,\end{array}\right.
	\end{equation}
	with boundary data $f$, given by the Riesz representation formula
	\begin{equation}\label{eq.riesz}
	u(X)=\int_{\partial\Omega}f(\xi)\,d\omega^X(\xi),\qquad X\in\Omega,
	\end{equation}
	satisfies the estimate 
	\begin{equation}\label{eq.direst}
		\Vert\m N(u)\Vert_{L^{p}(\partial\Omega,\sigma)}\leq C\,\Vert f\Vert_{L^{p}(\partial\Omega,\sigma)}.
	\end{equation}
	We call  $C$ in (\ref{eq.direst}) the \emph{$(\Di_{p}^L)$ constant}.  
	
	We say that the   \emph{Dirichlet regularity problem}   (or just \textit{regularity problem})  for the operator $L$ with   $\dt W^{1,p}(\bb R^n)$  data  is solvable   (write $(\Reg_{p}^L)$ is solvable), if there exists $C\geq1$ so that for each $f\in C_c^1(\bb R^n)$, the solution $u$ to the Dirichlet problem for $L$ (\ref{eq.dirichlet}) with boundary data $f$, given by the Riesz representation formula (\ref{eq.riesz}),  satisfies the estimate
	\begin{equation}\label{eq.regest}
		\Vert\wt{\m N}(\nabla u)\Vert_{L^{p}(\partial\Omega,\sigma)}\leq C\Vert\nabla f\Vert_{L^p(\bb R^n)}.
	\end{equation} 
	Here, $\dt W^{1,p}(\bb R^n)$ is the classical homogeneous Sobolev space of locally integrable functions in $\bb R^n$ whose weak gradient exists and lies in $L^p(\bb R^n)$. We call the constant $C$ in (\ref{eq.regest}) the \emph{$(\Reg_p^L)$ constant}.
\end{definition}

In \cite[Problem 3.2.25]{kbook}, Kenig asks whether for a given elliptic operator $L$ as defined above, one always has that 
\begin{equation}\label{eq.imply}
	(\Di_{p'}^{L})\implies(\Reg_p^L),
\end{equation}
where $p'$ is the H\"older conjugate of $p$, and $p\in(1,\infty)$. This has remained an outstanding open problem in the field, but there are several partial  and related results in the literature. We defer discussion of an important related result of Gallegos \cite{gallegos25} to Section \ref{sec.related}, and a brief historical review to Section \ref{sec.history}.  \\

We next state our main result, whereby we explicitly construct a real, symmetric (in fact, \emph{diagonal}), elliptic operator $L=-\dv A\nabla$ for which $(\Di_{p}^{L})$ is solvable for all $p\in(1,\infty)$, but $(\Reg_q^L)$ is not solvable for any $q>1$.



\begin{theorem}\label{main:counterexample}  There exists a real, diagonal, strongly elliptic matrix $A\in C^\infty(\Hn)\cap C(\overline{\Hn})$, with the following	properties.
	\begin{enumerate}
		\item $(\Di_p^L)$ is solvable for all $p\in(1,\infty)$, with the $(\Di_{p}^L)$ constant depending only on $n$ and $p$. 
		\item There exist $f\in C_c^\infty(\R^n)$ and a Borel set $F\subset\partial\Omega$ 
		with $0<\m H^n(F)<\infty$, such that
		\begin{equation}\label{intro:counter-failure}
			\widetilde {\m N}(\nabla u_f)(\xi)=\infty\qquad\text{for each }\xi\in F.
		\end{equation}
		In particular, $(\Reg_p^L)$ is not solvable for any $p>1$.
	\end{enumerate}
\end{theorem}

This example resolves Problem 3.2.25 of Kenig \cite{kbook}, and shows that in the full generality of real elliptic operators  on a smooth domain, $(\Di_{p'}^{L})$ does not imply $(\Reg_p^L)$.\\


Let us now state several interesting properties that the counterexample of Theorem \ref{main:counterexample} satisfies, and which are shown along the way as part of its construction. The properties highlighted below show that the counterexample lives quite close to conditions in the literature known to imply solvability of the Dirichlet regularity problem for at least some $p>1$. We defer all definitions to the preliminaries in Section \ref{sec:preliminaries}. 

\begin{theorem}\label{main2} Fix $\varepsilon\in(0,\frac14)$. There exists a real, diagonal, strongly elliptic matrix $A\in C^\infty(\Hn)\cap C(\overline{\Hn})$, satisfying the properties in Theorem \ref{main:counterexample}, and with the following additional properties.
\begin{enumerate}		
		\item We have that $A(x,0)$ equals the $(n+1)\times(n+1)$ identity matrix $I$; that is, $A|_{t=0}=I$. Moreover, $Ae_{n+1}=e_{n+1}$.
		
		\item We have that the scale-invariant local Lipschitz constant of $A$ is uniformly small and that $A$ is  uniformly close to the identity; namely,
		\begin{equation}\label{eq.smallperturb}
		\Vert A-I\Vert_{L^\infty(\bb R^{n+1}_+)}\leq\ep,\qquad \Vert t\nabla A\Vert_{L^\infty(\bb R^{n+1}_+)}\leq\ep.
		\end{equation}
		\item We have that  
		\begin{equation}\label{eq.vanishlip}
		\lim_{s\searrow0}~\sup_{x\in\bb R^n,\,0<t<s}t|\nabla A(x,t)|=0.
		\end{equation}
		\item The matrix satisfies the $t$-DKP condition (see Definition \ref{def.tdkp}), and in fact,
		\begin{multline}
			t\sup_{Y\in B((x,t), \frac t2)}
			|\partial_sA(Y)|^2\,dx\,dt\quad\text{is a Carleson measure}\\ \text{ with Carleson measure constant bounded above by }\ep^2.
			\label{intro:counter-gradient}
		\end{multline}
		
		\item The matrix is not a DKP matrix (see Definition \ref{def.dkp}), but for each $q>2$, we have that
		\begin{multline}\label{dr:vanishing-weak-quadratic}
			t^{q-1}\sup_{Y\in B((x,t),\frac t2)}|\nabla A(Y)|^q\,dx\,dt\quad\text{is a vanishing Carleson measure}\\ \text{ with Carleson measure constant bounded above by }C_q\ep^q.
		\end{multline}
		The constant  $C_q$ depends only on $n$ and $q$.
		
		
	\end{enumerate}
	
\end{theorem}

The proofs of Theorems \ref{main:counterexample} and \ref{main2} are shown in Section \ref{sec:counterexample-proof}. Let us compare the statements in Theorem \ref{main2} to corresponding results in the literature; see Section \ref{sec.history} for the historical review. The first estimate in (ii) tells us  that the counterexample matrix $A$ is an $L^\infty$ perturbation of the identity; but crucially, it is not $t$-independent, so it falls outside of the purview of the results in \cite{aaahk,bhlmp20,bhlmp24}. The second estimate in (ii) and the statement (iii) tell us that the scale-invariant local Lipschitz constant of $A$ is small, even vanishing, but that this is not enough to guarantee solvability of the regularity problem.\\

The statements (iv) and (v) both quantify, in different ways, how much the matrix $A$ fails to be a DKP matrix. Indeed, if $A$ were a DKP matrix, then it has recently been established that the regularity problem is solvable for some $p>1$ \cite{dhp22,mpt25}. Thus, it is not a surprise that $A$ fails to satisfy the DKP condition; however, statement (v) tells us that $A$ satisfies a DKP-type condition with all exponents larger than $2$, and thus $A$ barely fails to be a DKP matrix. Finally, statement (iv) gives that $A$ does satisfy a weaker DKP-type condition (with exponent $2$) whereby the Carleson measure condition is restricted to only the transversal derivative, and not the full gradient. Stronger conditions of this type have recently been considered in the literature by Ulmer \cite{ulmer2, ulmer3}, and so our result here shows that the implication (\ref{eq.imply}) cannot be true for $t$-DKP matrices, notwithstanding the results in \cite{dhp22,mpt25,ulmer2,ulmer3}.\\


In a similar spirit to results in the note \cite{pog1}, the construction of the counterexample moreover allows for an approximation of the matrix $A$ by a sequence of truncations $A_N\in C^\infty(\overline\Omega)$,   for each of which the regularity problem is solvable in all $p\in(1,\infty)$ (with non-uniform estimates, of course), and the Dirichlet problem is solvable with uniform estimates.  In this way, the counterexample lies at the boundary of the space of matrices for which the Regularity problem is solvable. These observations are recorded precisely in the following corollary.

\begin{corollary}\label{dr:approximations} Let $A$ be the matrix constructed in Theorem \ref{main:counterexample}. For each $N\in\bb N$, let $A_N$ be the truncation of the matrix retaining the frequencies $4^k$, $k\leq N$. For each $N\in\bb N$, $A_N$ is a DKP matrix (see Definition \ref{def.dkp}).  
The truncations $A_N$  converge
uniformly to $A$ and belong to $C^\infty(\overline{\Hn})$. If $L_N:=-\dv A_N\nabla$, then for each $p\in(1,\infty)$, the $(\Di_p^{L_N})$ constant is uniform in $N$. For each $N\in\bb N$,  $(\Reg_p^{L_N})$ is solvable for every $p\in(1,\infty)$. For the fixed
$f\in C_c^\infty(\R^n)$ in \eqref{intro:counter-failure}, and for every $p\in(1,\infty)$, we have that 
\[
 \lim_{N\to\infty}
 \|\widetilde{\m N}(\nabla u_{f,N})\|_{L^{p,\infty}(\R^n)}=\infty.
\] 
\end{corollary}

The proof of Theorem~\ref{main:counterexample} is inspired by a merging of the construction of related counterexamples for the Neumann problem in \cite{mourg26} (see also \cite{mpt26}) with a modification by the third author \cite{pog1} of the classical example of Modica-Mortola \cite{mm} regarding the existence of singular elliptic measures. \\

\subsection{Overview of the proof} We now give an overview of the proof of Theorem \ref{main:counterexample}.   In \cite{mourg26}, the author considers the following matrix: Working first in $\R^3_+$, with coordinates $(x_1,x_2,t)$ and $e_t=e_3$,  take $u=x_1+v(x_2,t)$. The skew coupling
\[
 A_{\mathrm{Neu}}\coloneqq I+b(x_2,t)(e_1\otimes e_t-e_t\otimes e_1)
\]
gives the equation and outward
conormal derivative
\[
 \Delta v=\partial_t b,\qquad
 \partial_{\nu_{A_{\mathrm{Neu}}}}u=(-v_t+b)|_{t=0}.
\]
In the local Neumann model, one imposes $v_t=b$ on $t=0$ and increases the tangential derivative of the boundary trace. The
Dirichlet datum therefore varies. Moreover, when $b$ is independent of $t$, the operator in the interior is the Laplacian, so this change of conormal derivative leaves the harmonic Dirichlet solutions and their regularity estimates unchanged. The full Neumann construction localizes this calculation and repeats an increase of the boundary-gradient integral, keeping one Neumann datum fixed. Here one uses the local reduction to a Poisson equation, but adds the corrections exactly instead of repeating that localization argument.

For our case at hand regarding the regularity problem, we instead impose $v(x_2,0)=0$, so the correction
does not change the boundary datum $x_1$. We also keep the normal
row and column equal to those of the identity, to preserve uniformity of the constant in the solvability of the Dirichlet problem $(\Di_p)$. We accordingly couple the first two
tangential directions and take 
\[
A=I+b(e_1\otimes e_2+e_2\otimes e_1)
\]
Since $b$ and $v$ are
independent of $x_1$, the equation for $u_0=x_1+v(x_2,t)$ becomes
\[
 v_{tt}+v_{x_2x_2}=-\partial_{x_2}b,\qquad v(x_2,0)=0.
\]

Our aim is to keep the coefficient perturbation small while
making the normal derivative of the solution large. To see how
a single oscillation contributes, take
\[
b(x_2,t)=a\chi(\lambda t)\cos(\lambda x_2),
\]
where $\chi\in C_c^\infty((1,2))$ is nonnegative and nonzero.
The corresponding correction has the form
\[
v(x_2,t)=\frac a\lambda\psi(\lambda t)\sin(\lambda x_2),
\qquad
\psi''-\psi=\chi,\quad \psi(0)=0,\quad
\lim_{s\to\infty}\psi(s)=0.
\]
Solving this one-dimensional equation gives
\[
\partial_tv(x_2,0)
=-ac_\chi\sin(\lambda x_2),
\qquad
c_\chi=\int_0^\infty e^{-s}\chi(s)\,ds>0.
\]
Thus increasing the frequency moves the coefficient layer
towards the boundary without reducing the size of its
contribution to the boundary normal derivative.

We will superpose such oscillations at sufficiently separated
frequencies. Their coefficient layers are disjoint, so the
size of the coefficient perturbation is controlled by the
largest amplitude. The solution corrections, however, extend
below those layers, and their normal derivatives accumulate.
At suitable positive observation heights, they reproduce
partial sums of a lacunary trigonometric series, up to a
uniformly bounded error. We therefore seek amplitudes tending
to zero whose squared sum diverges, so that the coefficients
approach the identity while these partial sums become
unbounded on a set of positive measure.

This use of oscillations in separated height layers is strongly related to the classical construction of Modica and Mortola \cite{mm} and its more recent modification by the third author \cite[Section~2 and Proposition~3.10]{pog1}. In the latter note, the author uses square-summable amplitudes and large horizontal frequencies to retain quantitative absolute continuity of the elliptic measure with respect to the surface measure, while the DKP constant degenerates to infinity. In that construction the
normal coefficient varies, and the Poisson kernels are controlled
by Riesz products. Here we keep the normal block fixed, which
gives all Dirichlet estimates independently of the oscillations.
The reduced equation is linear in $b$ and $v$. We can therefore
use amplitudes with divergent squared sum to obtain large normal
derivatives through an additive series.  

The fixed cutoff above, however, would give  
\[
 \int_0^\infty t\left|\partial_t\bigl(a\chi(\lambda t)\bigr)\right|^2
                  \dd t
 =a^2\int_0^\infty s|\chi'(s)|^2\dd s.
\]
Its vertical cost is comparable to $a^2$, so these costs cannot be
summed when $\sum_i a_i^2=\infty$. We instead spread the transitions
in the logarithm of the height. For a tangential frequency
$\lambda>0$, take
$b=a\zeta_L(\lambda t)\cos(\lambda x_2)$, where $\zeta_L$ is
zero outside $[e^{-2L},e^{2L}]$, is one on $[e^{-L},e^L]$,
and has transitions of length $L$ in the variable $\log t$.
The correction is
\[
 v(x_2,t)=\frac a\lambda\psi_L(\lambda t)\sin(\lambda x_2),
 \qquad
 \psi_L''-\psi_L=\zeta_L,\quad
 \psi_L(0)=0,\quad \lim_{s\to\infty}\psi_L(s)=0.
\]
The explicit solution gives
\[
 \partial_tv(x_2,0)=-ac_L\sin(\lambda x_2),\qquad
 c_L=\int_0^\infty e^{-s}\zeta_L(s)\dd s\geq\tfrac12
 \quad(L\geq2),
\]
whereas
\[
 \int_0^\infty s|\zeta_L'(s)|^2\dd s\leq C/L.
\]
Thus a long transition reduces the vertical derivative estimate
while retaining the size of the boundary normal derivative.

We superpose these corrections at separated frequencies which are
powers of $4$. The coefficient supports and their Whitney
enlargements are disjoint. The normal derivatives then contain
partial sums of the form $\sum_i a_ic_{L_i}\sin(\lambda_ix_2)$,
up to a uniformly bounded error at the selected heights.
The second and fourth moments give unbounded partial sums on a
set of positive measure as soon as $\sum_i a_i^2=\infty$.
Lemma~\ref{dr:boundary-growth} lets us replace $x_1$ by its product
with a smooth compactly supported cutoff. A fixed tangential
rotation makes the matrix diagonal.
Since $Ae_{n+1}=e_{n+1}$, comparison with the height function gives
all the Dirichlet estimates in Lemma~\ref{dr:all-dirichlet}.

It remains to choose the layers so that the coefficient estimates
hold simultaneously. We group them into finite blocks. In block
$j$, we use the logarithmic width $\mathcal L_j=L_*2^j$ and take
\[
M_j=\left\lceil e^{\mathcal L_j^2}\right\rceil
\]
layers, all with amplitude
\[
a_{j,\ell}=\frac{\delta}{\sqrt{M_j}},
\qquad 1\leq\ell\leq M_j.
\]
Here $L_*$ is a sufficiently large fixed dimensional constant,
and $\delta$ is a sufficiently small dimensional multiple of
$\varepsilon$. The frequencies are chosen recursively among
powers of $4$, sufficiently far apart that the height supports
of the layers and their Whitney enlargements are disjoint.

Each block has the same squared amplitude sum:
\[
\sum_{\ell=1}^{M_j}a_{j,\ell}^2=\delta^2.
\]
Consequently, the total squared amplitude sum diverges, as needed
for the unbounded normal derivatives, whereas the vertical
quadratic Carleson costs satisfy
\[
\sum_{j\geq1}\frac1{\mathcal L_j}
\sum_{\ell=1}^{M_j}a_{j,\ell}^2
=\delta^2\sum_{j\geq1}\frac1{\mathcal L_j}<\infty.
\]
The choices of $L_*$ and $\delta$ make the resulting vertical
Carleson norm at most $\varepsilon^2$.

The much larger number of layers in each successive block gives
the remaining estimates. The scaled full gradient contributed
by a layer in block $j$ is bounded by
\[
Ca_{j,\ell}e^{2\mathcal L_j}
\leq C\delta
\exp\left(2\mathcal L_j-\frac{\mathcal L_j^2}{2}\right),
\]
which is uniformly small and tends to zero with $j$. Moreover,
for each $q>2$, the contribution of block $j$ to the
$q$-power Carleson estimate is bounded by
\[
C_q\delta^q(1+\mathcal L_j)
\exp\left(
2q\mathcal L_j-\left(\frac q2-1\right)\mathcal L_j^2
\right).
\]
This is summable in $j$ for every $q>2$, with the same choice of layers. Since every finite collection of blocks is supported a positive distance from the boundary, the tails of these convergent series also give the vanishing Carleson property. 


\subsection{Related results}\label{sec.related} Gallegos \cite{gallegos25} recently proved  that, in the plane,  given a Corkscrew domain with Ahlfors regular boundary,  solvability of a \emph{generalized} regularity problem for the Laplacian in some $L^p$ implies uniform rectifiability of the boundary (see \cite{gallegos25} for the definitions of these geometric concepts). In addition, he showed that, for the example of David and Mayboroda \cite{dm2} of an elliptic operator for whom the Dirichlet problem is solvable in all $p\in(1,\infty)$ on a domain in the plane whose boundary has a purely unrectifiable part,  the generalized regularity problem is not solvable for any $p>1$.  The free boundary result of Gallegos has been extended to higher dimensions and general elliptic operators in divergence form by the first two authors and will appear in the manuscript  \cite{MMV}, using the resolution of the WALA conjecture by the first author (see \cite{MMV}). Moreover, the aforementioned example of David and Mayboroda was extended to higher dimensions by Bortz, Tapiola, Tolsa, and the third author \cite{bptt24}. Thus the implication (\ref{eq.imply}), properly interpreted in rough domains with the definition of a generalized regularity problem, does not hold for arbitrary elliptic operators in uniform domains with Ahlfors regular boundary, in any ambient dimension. \\

The free boundary result  in \cite{gallegos25} is very interesting and contributes deeply to the theory of the Dirichlet regularity problem in rough domains. Nevertheless, this  result  does not resolve Kenig's original question in \cite[Problem 3.2.25]{kbook}. This is because   the domains in which Kenig originally posed the problem of the implication (\ref{eq.imply}) were certainly uniformly rectifiable; and indeed, the problem was open even in the upper-half space and the unit ball. However, on smooth domains, or even more generally, on Corkscrew domains with uniformly rectifiable boundary, the argument in \cite{gallegos25} does not provide a counterexample to the implication (\ref{eq.imply}), because his argument relies heavily on the lack of rectifiability of the boundary.\\

\subsection{Historical results}\label{sec.history} On the upper-half space, the question of the solvability of the Dirichlet problem and the Dirichlet regularity problem  for general real  elliptic operators, as well as the relationship of solvability between both problems, has remained a topic of sustained study. It is classical that if $L$ is the Laplacian, then $(\Di_{p}^{-\Delta})$ and $(\Reg_p^{-\Delta})$ are solvable for all $p\in(1,\infty)$ (see for instance \cite{stein70}). If $L$ is symmetric and has $t$-independent coefficients,  then $(\Di_{2}^{L})$ and $(\Reg_2^{L})$ are solvable, by works of Jerison and Kenig, and Kenig and Pipher \cite{jk1, kp, kp2}. Moreover, it was  shown in \cite{kp} that for general elliptic operators, $(\Reg_p^L)\implies(\Di_{p'}^{L})$, but the converse direction remained an open problem.  Later, Shen proved in \cite{shen07} that on bounded Lipschitz domains, if $(\Di_{p'}^{L^*})$ is solvable and $(\Reg_q^L)$ is solvable for \emph{some} $q>1$, then $(\Reg_p^L)$ is solvable; and this result was subsequently extended by Dindo\v{s} and Kirsch \cite{dk12}. \\

The first non-trivial positive result for non-symmetric $t$-independent operators on the half-space $\bb R^{n+1}_+$ was obtained by Kenig, Koch, Pipher, and Toro \cite{kkpt} when $n=1$. There, the authors showed that there exists $p>1$ such that $(\Di_{p'}^L)$ is solvable, and unlike the symmetric $t$-independent case, $p$ may not be made precise.  Later, Kenig and Rule showed in \cite{kr} that $(\Di_{p'}^{L^*})\implies(\Reg_p^L)$ for the $t$-independent operators when $n=1$, and thus it followed that $(\Reg_p^L)$ was solvable for some $p>1$. These results of \cite{kkpt} and \cite{kr} were generalized to $n\geq2$ by Hofmann, Kenig, Mayboroda, and Pipher in \cite{hkmps} and \cite{hkmpr}, where the solution to the Kato problem \cite{hlm02, ahlmct} played a   role in the proof.      Let us also mention the work of Alfonseca, Auscher, Axelsson, Hofmann, and Kim \cite{aaahk}, in which the authors showed that the solvability of $(\Di_{2}^{L})$ and $(\Reg_2^{L})$ are preserved through  $t$-independent small  $L^\infty$  perturbations of the matrix at hand; this result would later be generalized to second-order elliptic operators with $t$-independent lower order terms in \cite{bhlmp20, bhlmp24}. \\

In view of the results cited above, the situation with $t$-independent real operators is relatively well-understood. However, for operators which do not satisfy $t$-independence, the picture is remarkably murkier. First, recall that in  the case of the Dirichlet problem with \emph{continuous} data,  one has solvability for all strongly elliptic $A$ with bounded measurable coefficients \cite{lsw}. In contrast, for the case of the Dirichlet problem with $L^p$ data, important counterexamples were obtained independently by Caffarelli, Fabes, and Kenig \cite{cfk} and Modica and Mortola \cite{mm}, of symmetric operators $L$ for which   $(\Di_{p}^{L})$ is not solvable for any $p>1$, whence the regularity problem is not solvable either for any $p>1$.\\

As mentioned in   \cite{kp3,dpp07},  Dahlberg posed two conjectures in 1984 regarding solvability of operators with coefficients possibly depending on the $t$ variable.   The first conjecture concerned whether the $L^p$-solvability of the Dirichlet problem was stable under certain perturbations of the coefficients. More precisely, if $L_0=-\dv A_0\nabla$ and $L=-\dv A\nabla$ are two  elliptic operators, such that $(\Di_p^{L_0})$ is solvable in the half-space $\bb R^{n+1}_+$ and
\begin{equation}\label{eq.fkp}
	\frac{\esssup_{y\in B(x,\delta(x)/2)}|A(y)-A_0(y)|^2}{\delta(x)}\,dx\quad\text{is a Carleson measure,}
\end{equation}
does it follow that $(\Di_q^L)$ is solvable in $\bb R^{n+1}_+$, for some $q>1$ possibly larger than $p$? This question was resolved in the affirmative by Fefferman, Kenig, and Pipher in \cite{fkp}. The condition (\ref{eq.fkp}) has come to be known as the \emph{Fefferman-Kenig-Pipher} condition (FKP).  The FKP condition is known to also preserve the solvability of $(\Reg_p^L)$ (with some $q$ possibly smaller than $p$) \cite{kp2}.\\


Dahlberg's second conjecture asked whether for some $p>1$, $(\Di_p^L)$ is solvable in the half-space when $L$ is a DKP operator (see Definition \ref{def.dkp}). This question was resolved in the affirmative by Kenig and Pipher in \cite{kp3}. Dindo\v{s}, Petermichl, and Pipher later proved in \cite{dpp07} that for any $p>1$, if $A$ is a (possibly non-symmetric) DKP matrix with small enough Carleson measure constant   and $\Omega$ is a Lipschitz domain with small enough Lipschitz constant, then $(\Di_p^L)$ is solvable in $\Omega$; under the same assumptions, $(\Reg_p^L)$ is also solvable in $\Omega$, but this was shown a decade later by Dindo\v{s}, Pipher, and Rule \cite{dpr17}. More recently, Dindo\v{s}, Hofmann, and Pipher \cite{dhp22}, and independently the second and third authors together with Tolsa \cite{mpt25}, showed the implication (\ref{eq.imply}) for DKP operators  in unbounded Lipschitz domains and in bounded Corkscrew domains with uniformly rectifiable boundary, respectively, thus establishing the solvability of the regularity problem for DKP operators in the upper-half space, in particular.\\


Finally, and back in the upper-half space, we mention the recent works of Ulmer \cite{ulmer26, ulmer1, ulmer2, ulmer3} which study solvability of boundary value problems for elliptic operators satisfying conditions that control the   transversal derivative of $A$. Our counterexample shows that a natural $L^2$ analogue of his condition in \cite{ulmer3}, which we call here the $t$-DKP condition (see Definition \ref{def.tdkp}), does not imply solvability of the regularity problem.\\


For further historical review, including solvability of boundary value problems for nonsymmetric operators, or complex operators in the upper-half space, or for the case of rough domains, please see \cite{mpt25, bhlmp24, dp24, gallegos25}.

\subsection{Outline} The rest of this manuscript is organized as follows. In Section \ref{sec:preliminaries} we introduce all relevant definitions and state useful technical lemmas. In Section \ref{sec:counterexample-proof} we give the proof  of  Theorem \ref{main:counterexample}.

\subsection{Statement on use of artificial intelligence.} Our starting point was the second author's construction of a counterexample for the Neumann problem for generic operators in the upper half-space \cite{mourg26}. In the local model there, the coefficients in positions $(1,n+1)$ and $(n+1,1)$ have opposite signs, so that the conormal derivative contains a tangential derivative of the solution. Here our idea was to couple the first two tangential directions by setting $A_{12}=A_{21}=b(x_2,t)$ and leaving all other entries equal to those of the identity. This choice preserves Dirichlet solvability and allows us to make the normal derivative large near the boundary for a fixed smooth Dirichlet datum. We used large language models to help with the calculations arising from this choice of coefficients, to work through details of the proofs, and to draft and revise the manuscript.   The authors take full responsibility of the mathematical content.

\section{Preliminaries}\label{sec:preliminaries}

Recall that a substantial portion of the preliminary setup has already been established at the beginning of Section \ref{sec:introduction}.

\subsection{DKP matrices and Carleson measures}
\label{pre:coefficient-conditions}

We write
\[
T(\xi,r)=B^n(\xi,r)\times(0,r),\qquad
dm(x,t)=\frac{dx\,dt}{t}.
\]
For a matrix smooth in the open half-space, we define
\begin{equation}\label{dr:nearness-D}
	D_A(x,t)
	=t\sup_{Y\in B((x,t),t/2)}
	|\nabla A(Y)|,
\end{equation}
and its discrepancy from the identity by
\begin{equation}\label{pre:discrepancy}
	E_A(x,t)
	=\sup_{Y\in B((x,t),t/2)}
	|A(Y)-I|.
\end{equation} 

\begin{definition}[Carleson measures and vanishing Carleson measures]\label{def.carleson} Given a measure $\mu$ on $\Omega$, we say that it is a \emph{Carleson measure} on $\Omega$ if
\[
\|\mu\|_{\m C}:=\sup_{\xi\in\partial\Omega,\ r>0}r^{-n}\mu(T(\xi,r))<\infty.
\]
We call $\Vert\mu\Vert_{\m C}$ its Carleson measure constant. Moreover, if we set
\[
\|\mu\|_{\m C(r_0)}:=\sup_{\xi\in\partial\Omega,~0<r\leq r_0}~r^{-n}\mu(T(\xi,r)),
\]
then we say that $\mu$ is a \emph{vanishing Carleson measure}, and write $\mu\in\m C_0$, if
\[
\lim_{r\searrow0}\Vert\mu\Vert_{\m C(r)}=0.
\]
\end{definition}

\begin{definition}[DKP matrix]\label{def.dkp} We say that a locally Lipschitz matrix $A$ satisfies the Dahlberg-Kenig-Pipher condition (DKP) if
\[
t\sup_{Y\in B((x,t), \frac t2)}
|\nabla A(Y)|^2\,dx\,dt\quad\text{is a Carleson measure}.
\]
\end{definition}

We shall also consider a weaker condition than the DKP condition, whereby the Carleson measure restriction is imposed only on the transversal derivative of the matrix. Similar conditions have been studied recently in \cite{ulmer26, ulmer1, ulmer2, ulmer3}.

\begin{definition}[$t$-DKP matrix]\label{def.tdkp} We say that a locally Lipschitz matrix $A$ satisfies the $t$-DKP condition if  
\[
t\sup_{Y\in B((x,t), \frac t2)}
|\partial_sA(Y)|^2\,dx\,dt\quad\text{is a Carleson measure}.
\]
Accordingly, we define
\begin{equation}\label{pre:vertical-carleson}
	\mathcal K_t(A)
	=\sup_{\xi\in\R^n,\,r>0}\frac1{r^n}
	\int_{T(\xi,r)}t\sup_{Y\in B((x,t), \frac t2)}
	|\partial_sA(Y)|^2\dd x\,\dd t.
\end{equation} 
\end{definition}


\subsection{Dirichlet solutions in the half-space}
\label{dr:prelim}

For $f\in C_c^\infty(\mathbb R^n)$, we choose a smooth compactly
supported extension $\widetilde f$ to $\overline{\Hn}$ and denote
by $\dot W^{1,2}_0(\Hn)$ the completion of
$C_c^\infty(\Hn)$ in the gradient norm. The Lax--Milgram theorem
gives a unique $w\in\dot W^{1,2}_0(\Hn)$ such that
\[
\int_{\Hn}A\nabla w\cdot\nabla\varphi
=-\int_{\Hn}A\nabla\widetilde f\cdot\nabla\varphi,
\qquad \varphi\in\dot W^{1,2}_0(\Hn).
\]
We set $u_f\coloneqq\widetilde f+w$. To check that this definition is
independent of the extension, we use the zero-trace difference of
two such solutions as an admissible test in the difference of
their weak equations. We obtain $\nabla(u_f-\widehat u_f)=0$,
and the solutions agree once the additive constant is fixed by
the zero trace or by the choice of homogeneous Sobolev
representative.

We solve the Dirichlet problem on $\Hn\cap B(0,R)$ with data $f$
on the flat part and zero on the hemisphere, where $R$ is large
enough that the support of $\widetilde f$ lies strictly inside
$B(0,R)$. We write $u_{f,R}=\widetilde f+w_R$ and note that the
correction has zero trace, since the extension has the prescribed
values on the whole boundary. Thus, we may test with $w_R$ and use
ellipticity to obtain
\[
\lambda\|\nabla w_R\|_2^2
\leq\frac1{\lambda}\|\nabla\widetilde f\|_2\|\nabla w_R\|_2,
\qquad
\|\nabla w_R\|_2
\leq(1/\lambda^2)\|\nabla\widetilde f\|_2.
\]
After extending $w_R$ by zero outside the half-ball, this estimate
bounds the extensions in $\dot W^{1,2}_0(\Hn)$ uniformly in $R$,
while the weak maximum principle gives
$|u_{f,R}|\leq\|f\|_\infty$. We may therefore pass by weak
compactness on bounded sets to a limit satisfying the weak
equation. By the boundary H\"older estimates on flat patches, this
limit has continuous boundary trace $f$, so uniqueness identifies
it with the variational solution above. In this way, we obtain
continuity up to the boundary and the maximum principle for that
solution, and the same approximation gives its elliptic-measure
representation.

For real symmetric matrices, we also use the following local
boundary comparison estimates with $\Omega=\Hn$. Here and below,
$\omega_\Omega$ and $G_\Omega$ denote elliptic measure and the
Green function for the matrix under consideration, normalized by
$L_{A,X}G_\Omega(X,Y)=\delta_Y$. Symmetry of the matrix gives
$G_\Omega(X,Y)=G_\Omega(Y,X)$. The Green function is positive and
satisfies
\begin{equation}\label{pre:green-bound}
	0\leq G_\Omega(X,Y)\leq C|X-Y|^{1-n};
\end{equation}
see \cite{gw}.
Let $Y\in\partial\Omega$ and $r>0$. If $u,v$ are positive solutions
in $\Omega\cap B(Y,8r)$, continuous and zero on
$\partial\Omega\cap B(Y,8r)$, then
\begin{equation}\label{dr:boundary-comparison}
	\frac{u(X)}{v(X)}
	\leq C\frac{u(a_{Y,r})}{v(a_{Y,r})},
	\qquad X\in\Omega\cap B(Y,r).
\end{equation}
Here $a_{Y,r}$ is a fixed interior point whose distances from
$Y$ and the boundary are comparable to $r$. We use this local
comparison only when the solutions have no pole in the larger
boundary ball. If a pole $X$ lies outside $B(Y,8r)$,
we also have
\begin{equation}\label{dr:measure-green}
	\omega_\Omega^X(\partial\Omega\cap B(Y,r))
	\leq Cr^{n-1}G_\Omega(a_{Y,r},X).
\end{equation}
The boundary comparison estimates are due to
Caffarelli--Fabes--Mortola--Salsa \cite{cfms}, and the Green
function bound is also proved in \cite{gw}. All these estimates
hold locally on flat boundary patches, with constants depending
only on dimension and ellipticity.



\subsection{Estimates for block matrices}
\label{dr:local}

\begin{lemma}\label{dr:all-dirichlet}
	Suppose $A$ is real symmetric and uniformly elliptic in $\Hn$,
	and
	\begin{equation}\label{dr:block}
		A(X)=\begin{pmatrix}B(X)&0\\0&1\end{pmatrix}.
	\end{equation}
	Then, $(D_p)_{L_A}$ is solvable for every $1<p<\infty$, with
	constants depending only on dimension, ellipticity, $p$, and aperture.
	Moreover, $t$ is a positive weak solution of $L_At=0$.
\end{lemma}

\begin{proof}
	Since $A\nabla t=e_t$, every compactly supported smooth test satisfies
	$\int A\nabla t\cdot\nabla\varphi=\int\partial_t\varphi=0$.
	Thus, the height function is a positive solution with zero
	boundary values.
	
	We fix a pole $X=(X',X_t)\in\Hn$ and a boundary point $(y,0)$,
	and we set $R=|X-(y,0)|$ and $r=R/64$. The function
	$G_{\Hn}(\cdot,X)$ has no pole in $B((y,0),8r)$.
	By boundary comparison with the height function, we obtain
	\begin{equation}\label{dr:first-boundary-factor}
		G_{\Hn}((y,s),X)
		\leq C\frac{s}{r}G_{\Hn}((y,r/2),X),
		\qquad 0<s<r/16.
	\end{equation}
	If $X_t\geq r/16$, then the distance between $(y,r/2)$ and $X$
	is at least $R-r/2$, and \eqref{pre:green-bound} gives
	\[
	G_{\Hn}((y,r/2),X)
	\leq CR^{1-n}
	\leq C\frac{X_t}{r}R^{1-n}.
	\]
	If $X_t<r/16$, we use symmetry and apply boundary comparison
	in the second variable at $(X',0)$. The pole $(y,r/2)$ is outside
	the larger comparison ball, since
	\[
	|(y,r/2)-(X',0)|
	\geq R-X_t-r/2>8r.
	\]
	Consequently,
	\[
	G_{\Hn}((y,r/2),X)
	\leq C\frac{X_t}{r}
	G_{\Hn}((y,r/2),(X',r/2))
	\leq C\frac{X_t}{r}R^{1-n}.
	\]
	In the last inequality, we used
	$|X'-y|\geq R-X_t\geq R-r/16$ and
	\eqref{pre:green-bound}. Substituting either estimate into
	\eqref{dr:first-boundary-factor} and using $r=R/64$, we obtain
	\[
	G_{\Hn}((y,s),X)
	\leq C\frac{sX_t}{R^{n+1}},
	\qquad 0<s<r/16.
	\]
	We now apply \eqref{dr:measure-green} at boundary radius $s$,
	with a reference point of height comparable to $s$. By the preceding
	Green function bound, after decreasing the fixed upper bound for
	$s/R$ if necessary, we obtain
	\begin{equation}\label{dr:measure-upper-small-ball}
		\omega_{\Hn}^X(B^n(y,s))
		\leq Cs^n\frac{X_t}{|X-(y,0)|^{n+1}},
		\qquad 0<s<cR.
	\end{equation}
	In particular, the right side is at most $CX_t^{-n}s^n$.
	We cover a Lebesgue-null set by balls of radii less than $cX_t$
	whose $n$th powers have arbitrarily small sum. Applying
	\eqref{dr:measure-upper-small-ball} to these balls and letting
	the sum tend to zero, we obtain
	$\omega_{\Hn}^X\ll\mathcal H^n|_{\mathbb R^n}$. By
	differentiation of measures, the density $k^X$ satisfies
	\begin{equation}\label{dr:halfspace-kernel}
		k^X(y)\leq C\frac{X_t}
		{(|X'-y|^2+X_t^2)^{(n+1)/2}}
		\quad\text{for almost every }y\in\mathbb R^n.
	\end{equation}
	
	Suppose that $X\in\gamma_{\Hn,a}(\xi)$, where we identify
	the boundary point $\xi$ with its tangential coordinates. Then,
	$|X'-\xi|\leq aX_t$, and hence
	\[
	X_t+|\xi-y|
	\leq(a+2)|X-(y,0)|.
	\]
	We split the boundary into $B^n(\xi,2X_t)$ and the annuli
	$2^jX_t<|y-\xi|\leq2^{j+1}X_t$, $j\geq1$.
	By \eqref{dr:halfspace-kernel}, the contribution of the $j$th
	annulus to the representation of $|u_f(X)|$ is at most
	\[
	C_a\frac{2^{-j}}{(2^jX_t)^n}
	\int_{B^n(\xi,2^{j+1}X_t)}|f(y)|\dd y
	\leq C_a2^{-j}\m Mf(\xi).
	\]
	The same bound with $j=0$ holds for the inner ball. Summing
	the estimates, we obtain $|u_f(X)|\leq C_aMf(\xi)$, uniformly
	on the cone. Hence, $\m Nu_f\leq C\m Mf$, and the
	Dirichlet estimates follow from the $L^p$ boundedness of the
	Hardy--Littlewood maximal operator. For smooth compactly supported
	data, the construction in Subsection~\ref{dr:prelim} identifies
	the elliptic-measure solution with the variational solution.
\end{proof}

The Dirichlet estimate above is also stated in
\cite[Theorem 2.5]{dp24}. Its hypotheses impose no derivative
or Carleson condition on the tangential block.

We will use the following lemma for the difference between the
solution with boundary datum $x_1$ and the solution with a compactly
supported datum.

\begin{lemma}\label{dr:boundary-growth}
	Suppose \eqref{dr:block} holds in $\Hn\cap B(Y,8r)$, where
	$Y\in\partial\Hn$. Let $h$ be a bounded real weak solution there,
	continuous and zero on the flat boundary, and assume $|h|\leq M_h$.
	Then
	\begin{equation}\label{dr:growth}
		|h(x,t)|
		\leq CM_h\frac{t}{r},
		\qquad (x,t)\in\Hn\cap B(Y,2r).
	\end{equation}
	If $X=(x,t)\in B(Y,3r/2)$ and $0<t<r/4$, then
	\begin{equation}\label{dr:error-gradient}
		\left(\fint_{B(X,t/2)}|\nabla h|^2\right)^{1/2}
		\leq CM_h/r.
	\end{equation}
	The constants depend only on dimension and ellipticity.
\end{lemma}

\begin{proof}
	We set $D\coloneqq\Hn\cap B(Y,8r)$ and choose continuous boundary
	data between zero and one, vanishing on $\partial D\cap B(Y,4r)$
	and equal to one on $\partial D\setminus B(Y,6r)$. We denote by
	$w$ the corresponding solution in $D$. By the maximum principle,
	$0\leq w\leq1$. We use $h=0$ on the flat part and
	$|h|\leq M_h=M_hw$ on the hemisphere to see that $h-M_hw$ and
	$-h-M_hw$ have nonpositive boundary values, and apply the
	maximum principle to these differences to obtain
	\[
	|h|\leq M_hw\quad\text{in }D.
	\]
	For $(x,t)\in B(Y,2r)$ with $t\geq r/8$, the bound $w\leq1$
	gives $w(x,t)\leq8t/r$. If $t<r/8$, we put $Y_x=(x,0)$
	and compare $w$ with the height function at $Y_x$, with reference
	radius $r/4$. The half-ball of radius $2r$ is contained in $D$,
	and its flat boundary lies in $B(Y,4r)$, where $w=0$, so boundary
	comparison gives
	\[
	w(x,t)
	\leq C\frac{t}{r}
	w\bigl(Y_x+(r/8)e_t\bigr)
	\leq C\frac{t}{r}.
	\]
	Combining the two bounds for $w$ near $Y$ with the preceding
	comparison for the signed differences, we get \eqref{dr:growth}.
	
	We choose $\phi\in C_c^\infty(B(X,3t/4))$ equal to one on
	$B(X,t/2)$, with $|\nabla\phi|\leq C/t$. We write $\lambda$ and
	$M$ for the lower ellipticity and upper matrix bounds. Since the
	support of the cutoff is contained in the domain, we may test the
	weak equation with $h\phi^2$ and, using ellipticity and Young's
	inequality, we get
	\[
	\begin{split}
		\lambda\int\phi^2|\nabla h|^2
		&\leq2M\int\phi|h||\nabla h||\nabla\phi|\\
		&\leq\frac\lambda2\int\phi^2|\nabla h|^2
		+\frac{2M^2}{\lambda}\int|h|^2|\nabla\phi|^2.
	\end{split}
	\]
	Since the height on $B(X,3t/4)$ is at most $7t/4$ and this ball
	is contained in $B(Y,2r)$, by \eqref{dr:growth} we can bound
	$|h|$ there by $CM_ht/r$. Absorbing the gradient term on the right
	of the preceding inequality and dividing by $|B(X,t/2)|$, we obtain
	\[
	\fint_{B(X,t/2)}|\nabla h|^2
	\leq\frac{C}{t^2}\fint_{B(X,3t/4)}|h|^2
	\leq CM_h^2/r^2.
	\]
\end{proof}

\section{Construction of the counterexample}\label{sec:counterexample-proof}

We prove Theorem \ref{main:counterexample}, Theorem \ref{main2}, and Corollary \ref{dr:approximations} with one construction.
Throughout the proof we use the fixed cones
\[
\gamma(q)=\{(x',t)\in\Hn:|(x'-q,t)|<2t\},
\]
and write $\m N$ and $\widetilde{\m N}$ for the corresponding maximal
functions; the averaging balls in $\widetilde{\m N}$ have radius $t/2$.

We first construct symmetric matrices coupling the first two
tangential directions. Their normal row and column agree with those
of the identity, so Lemma~\ref{dr:all-dirichlet} gives all Dirichlet
estimates. A fixed tangential rotation will make the matrices diagonal.

The coefficients are supported in separated horizontal layers.
We group the layers into finite blocks and give every layer within a
block the same amplitude and logarithmic width. Each block has squared
amplitude sum $\delta^2$. The widths increase so that the vertical
quadratic Carleson norm is finite. The number of layers in each block
increases much faster, which gives all the higher-power Carleson
estimates while the sum of squared amplitudes diverges.

\subsection{The cutoff and the one-frequency solution}
\label{uc:profiles}

Write $x=(x_1,x_2,z,t)$, where $z\in\R^{n-2}$ is omitted when
$n=2$. We require $v(x_2,0)=0$, so that $u_0=x_1+v(x_2,t)$ has
boundary value $x_1$, and keep the normal row and column of $A$
equal to those of the identity. Since $\partial_1u_0=1$, we choose
$Ae_1=e_1+b e_2$ to obtain a forcing term in the second tangential
direction. Taking the other entries from the identity and imposing
symmetry gives $Ae_2=e_2+b e_1$. Thus we set
\[
S=e_1\otimes e_2+e_2\otimes e_1,\qquad A=I+b(x_2,t)S.
\]
The first tangential variable supplies the constant derivative and
the second carries the oscillation. This is why we require $n\geq2$.
For smooth $b$ and $v$, direct multiplication gives
\begin{equation}\label{uc:generic-flux}
	A\nabla u_0=(1+bv_{x_2},\ b+v_{x_2},\ 0_z,\ v_t).
\end{equation}
Both $b$ and $v$ are independent of $x_1$ and $z$, so
\begin{equation}\label{uc:reduced-equation}
	\divg(A\nabla u_0)=b_{x_2}+v_{x_2x_2}+v_{tt}.
\end{equation}
Thus we solve $v_{tt}+v_{x_2x_2}=-b_{x_2}$ with $v(x_2,0)=0$.
The product $bv_{x_2}$ occurs only in the first flux component,
and its $x_1$ derivative is zero. The reduced equation therefore
contains no product of $b$ and $v$. This permits the exact
superposition below.

\DRAffineFigure

Fix a nondecreasing $\vartheta\in C^\infty(\R)$ with
$0\leq\vartheta\leq1$, $\vartheta=0$ on $(-\infty,0]$, and $\vartheta=1$
on $[1,\infty)$. For $L\geq2$, set
\begin{equation}\label{uc:log-cutoff}
	\zeta_L(s)
	=\vartheta\left(\frac{\log s+2L}{L}\right)
	\vartheta\left(\frac{2L-\log s}{L}\right)
	\quad(s>0).
\end{equation}
This function is zero for $s\leq e^{-2L}$ and for $s\geq e^{2L}$.
It is one on $[e^{-L},e^L]$. We extend it by zero at $s=0$.
This gives a smooth function on $[0,\infty)$.

The two transition intervals are disjoint. In each one, only one
factor in \eqref{uc:log-cutoff} is differentiated and the other
factor is one. Therefore,
\begin{equation}\label{uc:cutoff-bounds}
	|s\zeta_L'(s)|\leq C_\vartheta/L,\qquad
	\int_0^\infty s|\zeta_L'(s)|^2\dd s
	=\frac2L\int_0^1|\vartheta'(r)|^2\dd r.
\end{equation}
For the equality, use the variable $\log s$ on each transition
interval. The first inequality will also control the supremum
over Whitney balls.

For one coefficient term
$b(y,t)=a\zeta_L(\lambda t)\cos(\lambda y)$, with $y=x_2$, the
right-hand side of the reduced equation is
$a\lambda\zeta_L(\lambda t)\sin(\lambda y)$. We therefore choose
\begin{equation}\label{uc:single-mode}
	v(y,t)=\frac a\lambda\psi_L(\lambda t)\sin(\lambda y).
\end{equation}
Substitution gives
\[
v_{tt}+v_{yy}
=a\lambda\bigl(\psi_L''(\lambda t)-\psi_L(\lambda t)\bigr)
\sin(\lambda y).
\]
We consequently solve
\begin{equation}\label{uc:profile-equation}
	\psi_L''-\psi_L=\zeta_L,\qquad
	\psi_L(0)=0,\qquad \lim_{s\to\infty}\psi_L(s)=0.
\end{equation}
The solution is
\begin{equation}\label{uc:profile-formula}
	\psi_L(s)
	=-e^{-s}\int_0^s\sinh r\,\zeta_L(r)\dd r
	-\sinh s\int_s^\infty e^{-r}\zeta_L(r)\dd r.
\end{equation}
Indeed, differentiation gives
\[
\psi_L'(s)
=e^{-s}\int_0^s\sinh r\,\zeta_L(r)\dd r
-\cosh s\int_s^\infty e^{-r}\zeta_L(r)\dd r.
\]
The endpoint terms cancel in this derivative. Differentiating once
more gives $\psi_L''=\psi_L+\zeta_L$. Formula
\eqref{uc:profile-formula} also gives the two boundary conditions.
Uniqueness follows because a homogeneous solution vanishing at
zero is a multiple of $\sinh s$, and decay at infinity forces
that multiple to be zero.

We record bounds whose constants do not depend on $L$. The
nonnegative Green kernel for $-d^2/ds^2+1$ on the half-line is
\[
G(s,r)=\sinh(\min(s,r))e^{-\max(s,r)}.
\]
Since $0\leq\zeta_L\leq1$ and
$\int_0^\infty G(s,r)\dd r=1-e^{-s}$, we have
\begin{equation}\label{uc:profile-global}
	|\psi_L(s)|\leq1-e^{-s}\leq\min(s,1).
\end{equation}
The two nonnegative terms in the formula for $\psi_L'$ are
bounded respectively by
\[
e^{-s}(\cosh s-1)\leq\tfrac12,\qquad
\cosh s\,e^{-s}\leq1.
\]
Their difference therefore satisfies
\begin{equation}\label{uc:profile-derivative-global}
	|\psi_L'(s)|\leq1.
\end{equation}

Put
\[
c_L=\int_0^\infty e^{-r}\zeta_L(r)\dd r.
\]
For $L\geq2$,
\begin{equation}\label{uc:cL}
	1\geq c_L
	\geq e^{-e^{-L}}-e^{-e^L}\geq\tfrac12.
\end{equation}
The lower bound follows by integrating over $[e^{-L},e^L]$.
Below the support of $\zeta_L$, formula
\eqref{uc:profile-formula} becomes
\begin{equation}\label{uc:profile-low}
	\psi_L(s)=-c_L\sinh s,\qquad
	\psi_L'(s)=-c_L\cosh s
	\quad(0\leq s\leq e^{-2L}).
\end{equation}
Above the support, the same formula gives
\[
\psi_L(s)=-d_Le^{-s},\qquad
d_L=\int_0^{e^{2L}}\sinh r\,\zeta_L(r)\dd r
\leq e^{e^{2L}}.
\]
Consequently, for every integer $m\geq0$,
\begin{equation}\label{uc:profile-high}
	|\psi_L^{(m)}(s)|
	\leq \exp\bigl(-(s-e^{2L})\bigr)
	\qquad(s\geq e^{2L}).
\end{equation}
The derivatives in this range are all equal to $d_Le^{-s}$ up
to sign. No estimate for high derivatives in the transition
intervals is needed below.

For the function in \eqref{uc:single-mode}, the boundary normal
derivative is $-ac_L\sin(\lambda y)$. Thus a long logarithmic
transition reduces \eqref{uc:cutoff-bounds}, while \eqref{uc:cL}
keeps the boundary derivative comparable to $a$.

\subsection{Estimates for one coefficient layer}\label{uc:one-layer}

For $L\geq2$, $a\geq0$ and $\lambda\geq1$, consider the matrix
\[
P(x',t)=a\zeta_L(\lambda t)\cos(\lambda x_2)S.
\]
Write
\[
b_-=e^{-2L}/\lambda,\qquad b_+=e^{2L}/\lambda,
\qquad I^*=(2b_-/3,2b_+).
\]
This matrix is zero unless $b_-<t<b_+$. It extends smoothly by zero
across $t=0$, and
\begin{equation}\label{uc:layer-lipschitz}
	\|P\|_\infty\leq Ca,\qquad
	\|\nabla P\|_\infty\leq Ca/b_-.
\end{equation}
Indeed, a tangential derivative has size at most $Ca\lambda$, while
\eqref{uc:cutoff-bounds} bounds a vertical derivative by $Ca/(Lt)$.
On the support, $t\geq b_-$ and $\lambda\leq b_-^{-1}$.

A ball $B((x',t),t/2)$ can meet this support only when $t\in I^*$.
Every height $s$ in the ball satisfies $t/2<s<3t/2$. Consequently,
\begin{equation}\label{uc:layer-whitney}
	\begin{split}
		t\operatorname*{ess\,sup}_{(y',s)\in B((x',t),t/2)}
		|\nabla P(y',s)|
		&\leq Ca(\lambda t+L^{-1})\mathbf1_{I^*}(t)\\
		&\leq Ca e^{2L}\mathbf1_{I^*}(t),\\
		t\operatorname*{ess\,sup}_{(y',s)\in B((x',t),t/2)}
		|\partial_sP(y',s)|
		&\leq Ca L^{-1}\mathbf1_{I^*}(t).
	\end{split}
\end{equation}
The length of $I^*$ with respect to $dt/t$ is $4L+\log3$.
Integrating the second inequality gives
\begin{equation}\label{uc:layer-vertical-integral}
	\int_0^\infty t
	\left(\operatorname*{ess\,sup}_{(y',s)\in B((x',t),t/2)}
	|\partial_sP(y',s)|\right)^2\dd t
	\leq Ca^2/L.
\end{equation}
Squaring the first line of \eqref{uc:layer-whitney} and integrating
against $dt/t$ gives
\begin{equation}\label{uc:layer-quadratic-integral}
	\int_0^\infty
	\left(t\operatorname*{ess\,sup}_{(y',s)\in B((x',t),t/2)}
	|\nabla P(y',s)|\right)^2\frac{\dd t}{t}
	\leq C(1+L)e^{4L}a^2.
\end{equation}
For example, the two terms to integrate are bounded by
$Ca^2\lambda^2\int_0^{2b_+}t\,\dd t$ and
$Ca^2L^{-2}(4L+\log3)$.
All these estimates are uniform in $x'$.

Fix $C_0\geq1$ large enough to dominate the constants in
\eqref{uc:layer-lipschitz}--\eqref{uc:layer-quadratic-integral},
including the factor $|B^n(0,1)|$ from integration over boundary balls.
It depends only on dimension and the fixed smooth cutoff.

\subsection{The choice of the layers}\label{uc:selection}

Increasing $C_0$ if necessary to absorb the dimensional constants
in the coefficient estimates below, choose a fixed $L_*\geq4$ so large that
\[
C_0\sum_{j\geq1}(L_*2^j)^{-1}\leq1.
\]
These choices depend only on dimension and the fixed cutoff. Set
\[
\delta=\frac{\varepsilon}{10C_0},\qquad
\mathcal L_j=L_*2^j,\qquad
M_j=\left\lceil e^{\mathcal L_j^2}\right\rceil
\quad(j\geq1).
\]
Block $j$ consists of $M_j$ layers, all with logarithmic width
$\mathcal L_j$ and amplitude $\delta/\sqrt{M_j}$.

Enumerate the layers block by block with one index $i\geq1$, and let
$\mathcal I_j$ be the finite set of indices in block $j$. Thus
\begin{equation}\label{uc:block-properties}
	L_i=\mathcal L_j,\qquad a_i=\frac{\delta}{\sqrt{M_j}}
	\quad(i\in\mathcal I_j),\qquad
	\sum_{i\in\mathcal I_j}a_i^2=\delta^2.
\end{equation}
In particular,
\begin{equation}\label{uc:amplitudes-diverge}
	0<a_i\leq\delta,\qquad a_i\longrightarrow0,\qquad
	\sum_{i\geq1}a_i^2=\infty.
\end{equation}
The amplitude relevant to the scaled full gradient satisfies
\begin{equation}\label{uc:block-gradient-cap}
	a_i e^{2L_i}
	=\frac{\delta e^{2\mathcal L_j}}{\sqrt{M_j}}
	\leq\delta e^{2\mathcal L_j-\mathcal L_j^2/2}
	\leq\delta\qquad(i\in\mathcal I_j),
\end{equation}
and the bound tends to zero as $j\to\infty$.

We now choose the frequencies recursively. For positive integers
$m_i$, put
\begin{equation}\label{uc:active-notation}
	\lambda_i=4^{m_i},\qquad
	b_i=e^{-2L_i}/\lambda_i,\qquad
	T_i=e^{2L_i}/\lambda_i.
\end{equation}
Choose $m_1$ so that $T_1<1/4$. After choosing $m_i$, choose
$m_{i+1}>m_i$ so large that $T_{i+1}<b_i/100$. This is always
possible, and gives
\begin{equation}\label{uc:height-separation}
	T_1<1/4,\qquad T_{i+1}<b_i/100,\qquad
	T_i\leq T_1\,100^{-(i-1)}.
\end{equation}
Moreover,
\[
\frac{\lambda_{i+1}}{\lambda_i}
>100e^{2L_i+2L_{i+1}}>32.
\]
Consequently,
\begin{equation}\label{uc:frequency-separation}
	\lambda_{i+1}\geq32\lambda_i,\qquad
	\lambda_i>4\sum_{h<i}\lambda_h,
\end{equation}
since $\sum_{h<i}\lambda_h\leq\lambda_i/31$.

\subsection{The matrix and its coefficient estimates}\label{uc:coefficients}

Define
\begin{equation}\label{uc:matrix-construction}
	b(x_2,t)=\sum_{i\geq1}a_i\zeta_{L_i}(\lambda_it)
	\cos(\lambda_ix_2),\qquad A=I+bS.
\end{equation}
The height supports $(b_i,T_i)$ are disjoint, so at each point at
most one term is nonzero and $|b|\leq\delta$. By
\eqref{uc:height-separation}, the sum is locally finite in the open
half-space. Hence $A\in C^\infty(\Hn)$. Since $a_i\to0$, the
matrix tends uniformly to $I$ as $t\downarrow0$. Defining
$A(x',0)=I$ gives $A\in C(\overline{\Hn})$.

The matrix is symmetric, and
\begin{equation}\label{uc:ellipticity}
	Ae_{n+1}=e_{n+1},\qquad
	(1-\delta)|\xi|^2\leq A\xi\cdot\xi
	\leq(1+\delta)|\xi|^2.
\end{equation}
Indeed, $A\xi\cdot\xi=|\xi|^2+2b\xi_1\xi_2$ and
$2|\xi_1\xi_2|\leq|\xi|^2$. The choice of $\delta$ gives
$\|A-I\|_\infty\leq\varepsilon$ and uniform ellipticity.
Lemma~\ref{dr:all-dirichlet} gives $(D_p)_{L_A}$ for every
$1<p<\infty$, with constants depending only on dimension and $p$.

For each layer put $I_i^*=(2b_i/3,2T_i)$. These intervals are
disjoint by \eqref{uc:height-separation}. A Whitney ball can
therefore meet at most one coefficient layer. By
\eqref{uc:layer-whitney} and \eqref{uc:block-gradient-cap},
\begin{equation}\label{uc:scaled-whitney-cap}
	D_A(x',t)\leq C_0a_i e^{2L_i}
	\leq C_0\delta e^{2\mathcal L_j-\mathcal L_j^2/2}
	\quad(t\in I_i^*,\ i\in\mathcal I_j),
\end{equation}
and $D_A=0$ outside these intervals. Thus
$t|\nabla A|\leq D_A\leq C_0\delta\leq\varepsilon$.
Every finite collection of blocks has a positive bottom height,
and the last bound in \eqref{uc:scaled-whitney-cap} tends to zero
with $j$. It follows that
\begin{equation}\label{uc:scaled-gradient-vanishing}
	\lim_{r\downarrow0}\sup_{x'\in\R^n,\,0<t<r}
	t|\nabla A(x',t)|=0.
\end{equation}

\subsubsection{The vertical quadratic Carleson norm}

The estimate \eqref{uc:layer-vertical-integral}, disjointness of
the intervals $I_i^*$, and \eqref{uc:block-properties} give
\[
\begin{aligned}
	&\sup_{x'\in\R^n}\int_0^\infty t
	\left(\operatorname*{ess\,sup}_{(y',s)\in B((x',t),t/2)}
	|\partial_sA(y',s)|\right)^2\dd t\\
	&\hspace{2cm}\leq C\sum_{j\geq1}\frac1{\mathcal L_j}
	\sum_{i\in\mathcal I_j}a_i^2
	=C\delta^2\sum_{j\geq1}\mathcal L_j^{-1}.
\end{aligned}
\]
Integrating over $B^n(\xi,r)$, restricting the height integral to
$(0,r)$, and using the choices of $C_0,L_*$ gives
\begin{equation}\label{uc:vertical-smallness}
	\mathcal K_t(A)
	\leq C_0\delta^2\sum_{j\geq1}\mathcal L_j^{-1}
	\leq\delta^2\leq\varepsilon^2.
\end{equation}
This is the full vertical estimate with the supremum over Whitney
balls in the definition of $\mathcal K_t$.

\subsubsection{Higher-power Carleson estimates}

Fix $q>2$. Each interval $I_i^*$ in block $j$ has logarithmic
length $4\mathcal L_j+\log3$. Using
\eqref{uc:scaled-whitney-cap} and the equal amplitudes within
the block, we obtain, uniformly in $x'$,
\[
\begin{aligned}
	\int_{\bigcup_{i\in\mathcal I_j}I_i^*}
	D_A(x',t)^q\frac{\dd t}{t}
	&\leq C_q(1+\mathcal L_j)e^{2q\mathcal L_j}
	M_j\left(\frac{\delta}{\sqrt{M_j}}\right)^q\\
	&\leq C_q\delta^q(1+\mathcal L_j)
	\exp\left(2q\mathcal L_j
	-\left(\frac q2-1\right)\mathcal L_j^2\right).
\end{aligned}
\]
For each fixed $q>2$, the last expression is summable over $j$:
its exponent has a negative quadratic term in $\mathcal L_j$,
and $\mathcal L_j=L_*2^j$. Summing, integrating over boundary
balls, and absorbing the finite sum into $C_q$, we get
\begin{equation}\label{uc:higher-power-bound}
	\|D_A^qdm\|_{\mathcal C}
	\leq C_q\delta^q\leq C_q\varepsilon^q.
\end{equation}
The constants depend only on $q$, dimension, and the fixed cutoff.

These Carleson measures vanish on small scales. Indeed, for an
integer $J\geq1$, let
\[
\rho_J=\frac23\min_{i\in\bigcup_{j\leq J}\mathcal I_j}b_i>0.
\]
If a testing box has radius less than $\rho_J$, it misses the
enlarged intervals from the first $J$ blocks. Its normalized
$D_A^qdm$ measure is therefore bounded by the tail, over $j>J$,
of the convergent series above. This bound tends to zero as
$J\to\infty$. Hence $D_A^qdm\in\mathcal C_0$ for every $q>2$.

\subsubsection{Divergence of the quadratic norms}

Take $Q_*=(0,2\pi)^n$. Every layer lies below height $1/4$,
and $\zeta_{L_i}=1$ on $[1,2]$. Since $\lambda_i$ is an integer,
the averages of $\sin^2(\lambda_ix_2)$ and
$\cos^2(\lambda_ix_2)$ over $Q_*$ are both $1/2$.
The height supports are disjoint, so
\[
\begin{aligned}
	\int_{Q_*}\int_0^{1/4}t|\partial_{x_2}b|^2\dd t\dd x'
	&\geq\frac{|Q_*|}{2}\sum_i a_i^2\int_1^2s\,\dd s=\infty,\\
	\int_{Q_*}\int_0^{1/4}|b|^2\frac{\dd t\dd x'}{t}
	&\geq\frac{|Q_*|}{2}\sum_i a_i^2\int_1^2\frac{\dd s}{s}=\infty.
\end{aligned}
\]
Here we used \eqref{uc:amplitudes-diverge}. Since
$\nabla A=(\nabla b)S$ and $A-I=bS$, these integrals imply
divergence of the raw full-gradient and discrepancy integrals
in a fixed Carleson box containing $Q_*\times(0,1/4)$.
Continuity in the open half-space gives
$D_A\geq t|\nabla A|$ and $E_A\geq|A-I|$. Consequently,
\begin{equation}\label{uc:quadratic-discrepancy-infinite}
	\|t|\nabla A|^2\dd x'\dd t\|_{\mathcal C}
	=\|D_A^2dm\|_{\mathcal C}
	=\|E_A^2dm\|_{\mathcal C}=\infty.
\end{equation}

\subsection{An explicit solution}
\label{uc:solution}

Write $\psi_j=\psi_{L_j}$ and $c_j=c_{L_j}$. Define
\begin{equation}\label{uc:model-solution}
	v(y,t)=\sum_{j=1}^\infty\frac{a_j}{\lambda_j}
	\psi_j(\lambda_jt)\sin(\lambda_jy),
	\qquad
	u_0(x_1,x_2,z,t)=x_1+v(x_2,t).
\end{equation}
By \eqref{uc:profile-global} and
$\sum_j\lambda_j^{-1}<\infty$, the series for $v$ converges
uniformly on the closed half-plane. Its boundary value is zero.

We justify all interior differentiations. Fix $\tau>0$. For all
sufficiently large $j$, $T_j<\tau/2$. Thus, for $t\geq\tau$,
\[
\lambda_jt-e^{2L_j}\geq\lambda_j\tau/2.
\]
By \eqref{uc:profile-high}, each term with a total of $m$ derivatives
in $y,t$ is bounded by
\[
a_j\lambda_j^{m-1}e^{-\lambda_j\tau/2}.
\]
These bounds are summable because $\lambda_j\geq2^j$.
The finitely many remaining terms are smooth. Hence $v$ is
smooth for $t>0$, and its series can be differentiated there
to every order.

The profile equation \eqref{uc:profile-equation} and
\eqref{uc:single-mode} give
\[
v_{tt}+v_{yy}
=\sum_j a_j\lambda_j\zeta_{L_j}(\lambda_jt)
\sin(\lambda_jy)
=-b_y.
\]
The flux is given by \eqref{uc:generic-flux}.
Since $b$ and $v$ do not depend on $x_1$ or $z$, it follows
exactly that
\begin{equation}\label{uc:model-equation}
	\divg(A\nabla u_0)
	=\partial_{x_1}(1+bv_{x_2})
	+b_{x_2}+v_{x_2x_2}+v_{tt}=0.
\end{equation}

\subsubsection{Bounds near the boundary}

We prove, with constants independent of the chosen layers,
\begin{equation}\label{uc:model-bounds}
	\begin{split}
		|v_y(y,t)|&\leq C\delta,\\
		|v_t(y,t)|&\leq C\delta(1+\log_+(1/t)),\\
		|v(y,t)|&\leq C\delta t(1+\log_+(1/t)).
	\end{split}
\end{equation}
Here $\log_+s=\max(\log s,0)$.

For a fixed $t>0$, let $k$ be the largest index with $b_k\geq t$,
and set $k=0$ if there is no such index. For $i\leq k$,
$\lambda_it\leq e^{-2L_i}\leq1$. Equations
\eqref{uc:profile-global} and \eqref{uc:profile-derivative-global}
give
\[
\sum_{i\leq k}a_i|\psi_i(\lambda_it)|
\leq\delta t\sum_{i\leq k}\lambda_i
\leq C\delta t\lambda_k\leq C\delta,
\qquad
\sum_{i\leq k}a_i|\psi_i'(\lambda_it)|\leq\delta k.
\]
The sums are empty if $k=0$. For index $k+1$, both quantities
are at most $\delta$, by the same global profile bounds.

For $i\geq k+2$, the gap condition and $b_{k+1}<t$ imply
\[
T_i<\frac{t}{100^{\,i-k-1}},\qquad
\lambda_it>100^{\,i-k-1}e^{2L_i}.
\]
Thus \eqref{uc:profile-high} bounds both
$|\psi_i(\lambda_it)|$ and $|\psi_i'(\lambda_it)|$ by
\[
\exp\left(-(100^{\,i-k-1}-1)e^{2L_i}\right).
\]
Their sum over these indices is bounded by an absolute constant.
This proves the first bound in \eqref{uc:model-bounds} and
$|v_t|\leq C\delta(k+1)$.
If $k\geq1$, then
\[
t\leq b_k\leq T_k\leq T_1\,100^{-(k-1)}.
\]
Hence $k\leq C(1+\log_+(1/t))$, proving the second bound.
Finally, integrate this bound from $0$ to $t$ and use the
continuous boundary value $v(y,0)=0$. This gives the last bound.

It follows that $u_0$ belongs to $W^{1,2}$ on every bounded
cylinder up to $t=0$, because
$\int_0^1(1+\log(1/t))^2\dd t<\infty$. Its Sobolev trace
is $x_1$. The series in \eqref{uc:model-solution} also gives
$|v|\leq C\delta$ globally.

\subsubsection{The normal derivative at selected heights}

Set
\begin{equation}\label{uc:selected-heights}
	t_N=b_N/4,\qquad
	S_N(y)=\sum_{j=1}^N a_jc_j\sin(\lambda_jy).
\end{equation}
We claim that
\begin{equation}\label{uc:normal-partial-sums}
	|v_t(y,t)+S_N(\xi)|
	\leq C\delta
	\quad\text{if}\quad
	t_N/2\leq t\leq2t_N,\qquad |y-\xi|\leq t_N.
\end{equation}
For $j\leq N$, we have $t\leq b_N/2\leq b_j/2$.
Therefore \eqref{uc:profile-low} applies, and the $j$th
term of $v_t$ is
$-a_jc_j\cosh(\lambda_jt)\sin(\lambda_jy)$.
Its difference from $-a_jc_j\sin(\lambda_j\xi)$ is bounded by
\[
C a_j\left((\lambda_jt_N)^2+\lambda_jt_N\right).
\]
We used $c_j\leq1$, $|\sin r-\sin s|\leq|r-s|$ and
$|\cosh s-1|\leq Cs^2$ for $0\leq s\leq1$.
By \eqref{uc:frequency-separation},
$\sum_{j\leq N}\lambda_j\leq C\lambda_N$. Since
$\lambda_Nt_N=e^{-2L_N}/4\leq1/4$, these errors sum to at
most $C\delta$.

For $j>N$, iteration of \eqref{uc:height-separation} gives
$T_j<b_N100^{-(j-N)}$. Thus
\[
\lambda_jt\geq \frac{100^{\,j-N}}8e^{2L_j},\qquad
\lambda_jt-e^{2L_j}\geq 100^{\,j-N}/16.
\]
The high profile estimate gives
\[
\sum_{j>N}a_j|\psi_j'(\lambda_jt)|
\leq\delta\sum_{r\geq1}e^{-100^r/16}\leq C\delta.
\]
Combining the two estimates proves
\eqref{uc:normal-partial-sums}.

\DRObservationFigure

\subsection{Failure of the regularity estimates}
\label{uc:failure}

We first prove that the partial sums in
\eqref{uc:selected-heights} are unbounded on a set of positive
measure. Write $d_j=a_jc_j$ and
$s_N=\sum_{j\leq N}d_j^2$. By \eqref{uc:cL} and
\eqref{uc:amplitudes-diverge}, $s_N\to\infty$.
The integer frequencies are orthogonal on $(0,2\pi)$, and
\begin{equation}\label{uc:second-moment}
	\int_0^{2\pi}S_N(y)^2\dd y=\pi s_N.
\end{equation}
For the fourth moment we have exactly
\begin{equation}\label{uc:fourth-moment}
	\begin{split}
		\int_0^{2\pi}S_N(y)^4\dd y
		&=\frac{3\pi}{4}\sum_{j\leq N}d_j^4
		+3\pi\sum_{i<j\leq N}d_i^2d_j^2\\
		&\leq\frac{3\pi}{2}s_N^2.
	\end{split}
\end{equation}
To check the fourth-moment identity, expand each sine into two
exponentials. A nonzero integral requires a signed sum of four
frequencies to equal zero. By \eqref{uc:frequency-separation}, the
largest frequency cannot occur once, because it exceeds the sum
of three smaller frequencies. If it occurs twice with the same
sign, it cannot be cancelled by the other two terms. If it occurs
three times, their signed sum has magnitude at least that frequency,
which is larger than the remaining term. Thus only four equal
indices, or two equal pairs with opposite signs within each pair,
contribute. Their integrals are
\[
\int_0^{2\pi}\sin^4(\lambda_jy)\dd y=\frac{3\pi}{4},\qquad
\int_0^{2\pi}\sin^2(\lambda_iy)\sin^2(\lambda_jy)\dd y
=\frac\pi2\quad(i\ne j).
\]
The six orderings of each pair give the coefficient $3\pi$ in
\eqref{uc:fourth-moment}.

Let
$F_N=\{y\in(0,2\pi):|S_N(y)|\geq\frac12\sqrt{s_N}\}$.
The integral of $S_N^2$ over its complement is at most
$\pi s_N/2$. Equations \eqref{uc:second-moment} and
\eqref{uc:fourth-moment}, followed by Cauchy--Schwarz, give
\[
\frac{\pi s_N}{2}
\leq\int_{F_N}S_N^2
\leq |F_N|^{1/2}\left(\frac{3\pi}{2}\right)^{1/2}s_N.
\]
It follows that $|F_N|\geq\pi/6$. Since $s_N\to\infty$,
every set
\[
\{y\in(0,2\pi):\sup_N|S_N(y)|>m\},\qquad m\in\N,
\]
has measure at least $\pi/6$. These sets decrease with $m$.
Continuity of measure from above therefore shows that
\begin{equation}\label{uc:bad-partial-sums}
	F_*:=\{y\in(0,2\pi):\sup_N|S_N(y)|=\infty\}
	\quad\text{satisfies}\quad |F_*|\geq\pi/6.
\end{equation}

We now replace the boundary datum $x_1$ by a compactly supported
one. Put
\[
Q=(-1,1)\times(0,2\pi)\times(-1,1)^{n-2},
\]
with the last factor omitted when $n=2$. Choose
$\eta_0\in C_c^\infty(\R^n)$ equal to one on a fixed neighborhood
of $\overline Q$, and set
\begin{equation}\label{uc:compact-datum}
	f(x')=x_1\eta_0(x').
\end{equation}
Choose $r_*>0$ so that $\eta_0=1$ at all boundary points within
distance $9r_*$ of $\overline Q$.
Let $u_f$ be the Dirichlet solution, and set $h=u_f-u_0$.
The maximum principle gives $|u_f|\leq\|f\|_\infty$.
The model bounds \eqref{uc:model-bounds} bound $u_0$ on all
half-balls of radius $8r_*$ centered on $\overline Q$, uniformly
in the coefficient parameters. Thus $h$ is uniformly bounded
there, solves $Lh=0$, and has zero boundary values on their
flat parts. Lemma~\ref{dr:boundary-growth}, with $r=r_*$, gives,
for all sufficiently large $N$,
\begin{equation}\label{uc:error-averages}
	\left(\avint_{B^{n+1}((q,t_N),t_N/2)}
	|\nabla h|^2\right)^{1/2}\leq C_f
	\qquad(q\in Q).
\end{equation}
The constant is independent of $N$ and $q$.

Put $B_N(q)=B^{n+1}((q,t_N),t_N/2)$. Its center belongs to
$\gamma(q)$, and in fact the whole ball lies in $\gamma(q)$:
if $(x',t)\in B_N(q)$, then $|x'-q|<t_N/2<t$, so
$|(x'-q,t)|<\sqrt2\,t<2t$. Thus both maximal functions
dominate the same average. By \eqref{uc:normal-partial-sums}
and \eqref{uc:error-averages}, for all sufficiently large $N$,
\begin{equation}\label{uc:maximal-lower}
	\begin{aligned}
		\min\{\m N(\nabla u_f)(q),\widetilde{\m N}(\nabla u_f)(q)\}
		&\geq\left(\avint_{B_N(q)}|\nabla u_f|^2\right)^{1/2}\\
		&\geq\left|\avint_{B_N(q)}\partial_tu_f\right|
		\geq |S_N(q_2)|-C_f,\qquad q\in Q.
	\end{aligned}
\end{equation}

Define
\[
E_0=(-1,1)\times F_*\times(-1,1)^{n-2}\subset Q.
\]
Its measure is at least $2^{n-1}\pi/6$ by
\eqref{uc:bad-partial-sums}. At each point of $E_0$, the
partial sums remain unbounded after deleting any finite number
of indices. The bound \eqref{uc:maximal-lower}, which holds for all sufficiently
large $N$, proves
\[
\widetilde{\m N}(\nabla u_f)(q)=\m N(\nabla u_f)(q)=\infty
\quad(q\in E_0).
\]
For every $m>0$ and $1<p<\infty$,
\[
\|\m N(\nabla u_f)\|_{L^{p,\infty}}
\geq m|E_0|^{1/p},
\]
and the same inequality holds for $\widetilde N$. Letting
$m\to\infty$ rules out weak estimates, and hence strong ones.
The right hand side of the regularity estimate is finite, since
$\|f\|_{\dot W^{1,p}(\R^n)}=\|\nabla_Tf\|_{L^p(\R^n)}<\infty$.

\subsection{Smooth truncations and their regularity constants}
\label{uc:truncations}

For each natural number $N$, retain precisely the active frequency
exponents $m_i\leq N$ and put
\begin{equation}\label{uc:natural-truncations}
	b^{(N)}(x_2,t)=\sum_{m_i\leq N}a_i\zeta_{L_i}(\lambda_it)
	\cos(\lambda_ix_2),
	\qquad A_N=I+b^{(N)}S.
\end{equation}
At the unused exponents the coefficient is zero. Thus $N$ is the
natural exponent of the frequency $4^N$, and not a new numbering
of the nonzero terms. Every sum in \eqref{uc:natural-truncations}
is finite. It vanishes near $t=0$, so $A_N$ is smooth on the
closed half-space. Ellipticity, the normal block, and the
scaled-gradient bounds are
uniform in $N$. Lemma~\ref{dr:all-dirichlet} therefore gives all
Dirichlet estimates with constants independent of $N$.
Moreover,
\begin{equation}\label{uc:uniform-approximation}
	\|A_N-A\|_\infty\leq C\sup_{m_i>N}a_i\longrightarrow0.
\end{equation}
Here every block is finite and its amplitude $\delta/\sqrt{M_j}$
tends to zero as $j\to\infty$.

The truncations have a common modulus of continuity on the closed
half-space. To see this, let $P_i$ be the $i$th coefficient layer
in \eqref{uc:matrix-construction}. For each fixed $M$, put
$B_M=\sum_{i\leq M}\|\nabla P_i\|_\infty<\infty$, using
\eqref{uc:layer-lipschitz}. Since the height supports are disjoint,
for every $N$ and all $x,y\in\overline{\Hn}$ we have
\[
|A_N(x)-A_N(y)|
\leq B_M|x-y|+2C\sup_{i>M}a_i.
\]
First choosing $M$ large and then $|x-y|$ small proves uniform
equicontinuity. The same bound holds for $A$.

Every truncated full DKP norm is finite. Indeed, disjointness
of the enlarged height intervals and
\eqref{uc:layer-quadratic-integral} give
\begin{equation}\label{uc:finite-truncated-DKP}
	\|D_{A_N}^2dm\|_{\mathcal C}
	\leq C_0\sum_{m_i\leq N}(1+L_i)e^{4L_i}a_i^2<\infty.
\end{equation}
The sum is finite because only finitely many active layers are
retained. On the other hand, the same integral over
$Q_*\times(0,1/4)$ used in
\eqref{uc:quadratic-discrepancy-infinite} gives
\[
\|D_{A_N}^2dm\|_{\mathcal C}
\geq c\sum_{m_i\leq N}a_i^2\longrightarrow\infty.
\]
Every completed block contributes $\delta^2$ to the last sum.
This proves the desired behavior of the DKP condition of the truncations in Corollary \ref{dr:approximations}.

For each finite $N$, the normal block and the finite full
Whitney-supremum gradient norm give the averaged regularity estimate
for every $1<p<\infty$, by \cite[Theorem~5.6]{dp24}.
We prove that their constants must diverge using the same datum
$f$ in \eqref{uc:compact-datum}.

Let $r(N)$ be the number of retained active layers, and let
$v_N$ and $u_{0,N}=x_1+v_N(x_2,t)$ be obtained by retaining these
layers in \eqref{uc:model-solution}. The profile estimates and
\eqref{uc:height-separation} give all bounds in
\eqref{uc:model-bounds} uniformly for $v_N$: the omitted terms
only remove terms from the absolute-value estimates in their proof.
Let $u_{f,N}$ be the Dirichlet solution for $A_N$. On the fixed
boundary patches used above, $h_N=u_{f,N}-u_{0,N}$ has zero trace
and is bounded uniformly in $N$. Indeed, use the maximum principle
for $u_{f,N}$ and the uniform bound for $u_{0,N}$.
Lemma~\ref{dr:boundary-growth} therefore gives
\eqref{uc:error-averages} for $h_N$ with the same fixed constant.

For $r=r(N)$ sufficiently large, use the observation height
$t_r=b_r/4$. The proof of \eqref{uc:normal-partial-sums} applies
to the truncated model: all terms with index greater than $r$
are now absent, while the estimate for the first $r$ terms is
unchanged. Consequently, for every $\xi\in Q$,
\[
\widetilde{\m N}(\nabla u_{f,N})(\xi)
\geq |S_r(\xi_2)|-C,
\qquad
\m N(\nabla u_{f,N})(\xi)
\geq |S_r(\xi_2)|-C.
\]
By the second- and fourth-moment estimates, $|S_r|$ is at least
$\frac12\sqrt{s_r}$ on a set of measure at least $\pi/6$ in
$(0,2\pi)$. Multiplying by the other coordinate intervals of $Q$
gives, for every $1<p<\infty$,
\begin{equation}\label{uc:truncated-regularity-lower}
	\begin{split}
		\|\widetilde{\m N}(\nabla u_{f,N})\|_{L^{p,\infty}(\R^n)}
		&\geq c_{n,p}\sqrt{s_{r(N)}}-C_{n,p},\\
		\|\m N(\nabla u_{f,N})\|_{L^{p,\infty}(\R^n)}
		&\geq c_{n,p}\sqrt{s_{r(N)}}-C_{n,p}.
	\end{split}
\end{equation}
If a right-hand side is negative, the inequality is automatic;
thus increasing the constants covers the finitely many initial
indices as well. Since $c_i\geq1/2$ and every block has squared
mass $\delta^2$, we have $s_{r(N)}\to\infty$.
The same smooth compactly supported datum therefore gives the
claimed divergence for every exponent.
Its boundary-gradient norm is finite and fixed, so both the weak
and strong regularity constants tend to infinity.

\subsection{Diagonal form and completion of the proof}
\label{uc:diagonal}

Let $T$ be the orthogonal transformation defined by
\[
Te_1=(e_1+e_2)/\sqrt2,\qquad
Te_2=(-e_1+e_2)/\sqrt2,\qquad Te_j=e_j\quad(j\geq3).
\]
It preserves the half-space and its height coordinate. Put
\begin{equation}\label{uc:rotation}
	\widehat A(x)=T^{\mathsf T}A(Tx)T,\qquad
	\widehat A_N(x)=T^{\mathsf T}A_N(Tx)T,\qquad
	\widehat f(x')=f(Tx').
\end{equation}
The restriction of $T$ to the first $n$ coordinates is understood
in the last formula. Since $T^{\mathsf T}ST$ has diagonal entries
$1,-1,0,\ldots,0$, we obtain
\[
\widehat A(x)=\operatorname{diag}
(1+b(Tx),1-b(Tx),1,\ldots,1),
\]
and the same formula with $b^{(N)}$ gives $\widehat A_N$.
The normal entry remains $1$.

Changing variables $y=Tx$ in the weak equation shows that
$\widehat u_f=u_f\circ T$ solves
$-\divg(\widehat A\nabla\widehat u_f)=0$ with boundary datum
$\widehat f$. Uniqueness identifies it with the Dirichlet solution.
The same argument applies to every truncated matrix.
The transformation preserves balls, cones, height, Lebesgue measure
and gradient lengths, and
\[
\partial_t\widehat A(x)=T^{\mathsf T}(\partial_tA)(Tx)T,
\qquad |\nabla\widehat A(x)|=|\nabla A(Tx)|.
\]
It therefore preserves the common modulus of continuity, every
coefficient estimate, the vertical norm, the higher-power vanishing
Carleson properties, and the full quadratic and discrepancy divergences. It also preserves the natural truncations
and the finiteness and divergence of their DKP norms.
For the solutions,
\[
\m N(\nabla\widehat u_f)(\xi)=\m N(\nabla u_f)(T\xi),\qquad
\widetilde{\m N}(\nabla\widehat u_f)(\xi)
=\widetilde{\m N}(\nabla u_f)(T\xi),
\]
and similarly for the truncated solutions.
The set $F_*$ in \eqref{uc:bad-partial-sums} is Borel, since it is
defined by a countable supremum of continuous functions. Thus
$F=T^{-1}E_0$ is a Borel set of finite positive measure. Replacing $A,A_N,f$ by
$\widehat A,\widehat A_N,\widehat f$ proves
Theorem~\ref{main:counterexample} and Corollary~\ref{dr:approximations}.

\bibliographystyle{alpha-sort-max} 
\bibliography{refs} 

\begin{thebibliography}{AAAHK11}

\bibitem[AAAHK11]{aaahk}
M.~Alfonseca, P.~Auscher, A.~Axelsson, S.~Hofmann, and S.~Kim.
\newblock Analyticity of layer potentials and {$L^2$} solvability of boundary
  value problems for divergence form elliptic equations with complex
  {$L^\infty$} coefficients.
\newblock {\em Adv. Math.}, 226(5):4533--4606, 2011.

\bibitem[AHLMT02]{ahlmct}
P.~Auscher, S.~Hofmann, M.~Lacey, A.~McIntosh, and P.~Tchamitchian.
\newblock The solution of the {K}ato square root problem for second order
  elliptic operators on {${\mathbb R}^n$}.
\newblock {\em Ann. of Math. (2)}, 156(2):633--654, 2002.

\bibitem[BPTT24]{bptt24}
S.~Bortz, B.~Poggi, O.~Tapiola, and X.~Tolsa.
\newblock The {$A_\infty$} condition, {$\varepsilon$}-approximators, and
  {V}aropoulos extensions in uniform domains.
\newblock {\em J. Geom. Anal.}, 34(7):Paper No. 218, 53, 2024.

\bibitem[BHLMP22]{bhlmp20}
S.~Bortz, S.~Hofmann, J.~L. {Luna Garc\'ia}, S.~Mayboroda, and B.~Poggi.
\newblock Critical perturbations for second-order elliptic operators, {I}:
  square function bounds for layer potentials.
\newblock {\em Anal. PDE}, 15(5):1215--1286, 2022.

\bibitem[BHLMP24]{bhlmp24}
S.~Bortz, S.~Hofmann, J.~L. {Luna Garcia}, S.~Mayboroda, and B.~Poggi.
\newblock Critical perturbations for second order elliptic operators---{P}art
  {II}: {N}on-tangential maximal function estimates.
\newblock {\em Arch. Ration. Mech. Anal.}, 248(3):Paper No. 31, 96, 2024.

\bibitem[CFMS81]{cfms}
L.~Caffarelli, E.~Fabes, S.~Mortola, and S.~Salsa.
\newblock Boundary behavior of nonnegative solutions of elliptic operators in
  divergence form.
\newblock {\em Indiana Univ. Math. J.}, 30(4):621--640, 1981.

\bibitem[CFK81]{cfk}
L.~A. Caffarelli, E.~B. Fabes, and C.~E. Kenig.
\newblock Completely singular elliptic-harmonic measures.
\newblock {\em Indiana Univ. Math. J.}, 30(6):917--924, 1981.

\bibitem[DM21]{dm2}
G.~David and S.~Mayboroda.
\newblock Good elliptic operators on {C}antor sets.
\newblock {\em Adv. Math.}, 383:Paper No. 107687, 21, 2021.

\bibitem[DP24]{dp24}
M.~Dindo\v{s} and J.~Pipher.
\newblock Boundary value problems for elliptic operators satisfying {Carleson}
  condition.
\newblock {\em Vietnam J. Math.}, 52(3):627--673, 2024.

\bibitem[DHP23]{dhp22}
M.~Dindo\v{s}, S.~Hofmann, and J.~Pipher.
\newblock Regularity and {N}eumann problems for operators with real
  coefficients satisfying {C}arleson conditions.
\newblock {\em J. Funct. Anal.}, 285(6):Paper No. 110024, 32, 2023.

\bibitem[DK12]{dk12}
M.~Dindo\v{s} and J.~Kirsch.
\newblock The regularity problem for elliptic operators with boundary data in
  {H}ardy-{S}obolev space {$HS^1$}.
\newblock {\em Math. Res. Lett.}, 19(3):699--717, 2012.

\bibitem[DPP07]{dpp07}
M.~Dindo\v{s}, S.~Petermichl, and J.~Pipher.
\newblock The {$L^p$} {D}irichlet problem for second order elliptic operators
  and a {$p$}-adapted square function.
\newblock {\em J. Funct. Anal.}, 249(2):372--392, 2007.

\bibitem[DPR17]{dpr17}
M.~Dindo\v{s}, J.~Pipher, and D.~Rule.
\newblock Boundary value problems for second-order elliptic operators
  satisfying a {C}arleson condition.
\newblock {\em Comm. Pure Appl. Math.}, 70(7):1316--1365, 2017.

\bibitem[FKP91]{fkp}
R.~A. Fefferman, C.~E. Kenig, and J.~Pipher.
\newblock The theory of weights and the {D}irichlet problem for elliptic
  equations.
\newblock {\em Ann. of Math. (2)}, 134(1):65--124, 1991.

\bibitem[Gal]{gallegos25}
J.~M. Gallegos.
\newblock One-sided rellich inequalities, regularity problem and uniform
  rectifiability.
\newblock Preprint. 2025. arXiv: 2506.03431.

\bibitem[GW82]{gw}
M.~Gr\"uter and K.-O. Widman.
\newblock The {G}reen function for uniformly elliptic equations.
\newblock {\em Manuscripta Math.}, 37(3):303--342, 1982.

\bibitem[Ha96]{haj96}
P.~Haj\l~asz.
\newblock Sobolev spaces on an arbitrary metric space.
\newblock {\em Potential Anal.}, 5(4):403--415, 1996.

\bibitem[HKMP15a]{hkmpr}
S.~Hofmann, C.~Kenig, S.~Mayboroda, and J.~Pipher.
\newblock The regularity problem for second order elliptic operators with
  complex-valued bounded measurable coefficients.
\newblock {\em Math. Ann.}, 361(3-4):863--907, 2015.

\bibitem[HKMP15b]{hkmps}
S.~Hofmann, C.~Kenig, S.~Mayboroda, and J.~Pipher.
\newblock Square function/non-tangential maximal function estimates and the
  {D}irichlet problem for non-symmetric elliptic operators.
\newblock {\em J. Amer. Math. Soc.}, 28(2):483--529, 2015.

\bibitem[HLM02]{hlm02}
S.~Hofmann, M.~Lacey, and A.~McIntosh.
\newblock The solution of the {K}ato problem for divergence form elliptic
  operators with {G}aussian heat kernel bounds.
\newblock {\em Ann. of Math. (2)}, 156(2):623--631, 2002.

\bibitem[JK81]{jk1}
D.~S. Jerison and C.~E. Kenig.
\newblock The {D}irichlet problem in nonsmooth domains.
\newblock {\em Ann. of Math. (2)}, 113(2):367--382, 1981.

\bibitem[KKPT00]{kkpt}
C.~Kenig, H.~Koch, J.~Pipher, and T.~Toro.
\newblock A new approach to absolute continuity of elliptic measure, with
  applications to non-symmetric equations.
\newblock {\em Adv. Math.}, 153(2):231--298, 2000.

\bibitem[Ken94]{kbook}
C.~E. Kenig.
\newblock Harmonic analysis techniques for second order elliptic boundary value
  problems, volume~83 of {\em CBMS Regional Conference Series in Mathematics}.
\newblock Published for the Conference Board of the Mathematical Sciences,
  Washington, DC; by the American Mathematical Society, Providence, RI, 1994.

\bibitem[KP93]{kp}
C.~E. Kenig and J.~Pipher.
\newblock The {N}eumann problem for elliptic equations with nonsmooth
  coefficients.
\newblock {\em Invent. Math.}, 113(3):447--509, 1993.

\bibitem[KP95]{kp2}
C.~E. Kenig and J.~Pipher.
\newblock The {N}eumann problem for elliptic equations with nonsmooth
  coefficients. {II}.
\newblock {\em Duke Math. J.}, 81(1):227--250 (1996), 1995.
\newblock A celebration of John F. Nash, Jr.

\bibitem[KP01]{kp3}
C.~E. Kenig and J.~Pipher.
\newblock The {D}irichlet problem for elliptic equations with drift terms.
\newblock {\em Publ. Mat.}, 45(1):199--217, 2001.

\bibitem[KR09]{kr}
C.~E. Kenig and D.~J. Rule.
\newblock The regularity and {N}eumann problem for non-symmetric elliptic
  operators.
\newblock {\em Trans. Amer. Math. Soc.}, 361(1):125--160, 2009.

\bibitem[LSW63]{lsw}
W.~Littman, G.~Stampacchia, and H.~F. Weinberger.
\newblock Regular points for elliptic equations with discontinuous
  coefficients.
\newblock {\em Ann. Scuola Norm. Sup. Pisa Cl. Sci. (3)}, 17:43--77, 1963.

\bibitem[MMV]{MMV}
A.~Merlo, M.~Mourgoglou, and M.~Villa.
\newblock On the {WALA} conjecture, {A}lberti representations and applications.
\newblock Preprint. 2026. arXiv: 2605.29804.

\bibitem[MM81]{mm}
L.~Modica and S.~Mortola.
\newblock Construction of a singular elliptic-harmonic measure.
\newblock {\em Manuscripta Math.}, 33(1):81--98, 1980/81.

\bibitem[Mou]{mourg26}
M.~Mourgoglou.
\newblock On the ${L}^p$-{N}eumann problem for divergence form elliptic
  equations in the upper half-space.
\newblock Preprint. 2026. arXiv: 2611.24808.

\bibitem[MPT25]{mpt25}
M.~Mourgoglou, B.~Poggi, and X.~Tolsa.
\newblock Solvability of the {P}oisson-{D}irichlet problem with interior data
  in ${L}^{p'}$-{C}arleson spaces and its applications to the
  ${L}^p$-regularity problem.
\newblock {\em J. Eur. Math. Soc.}, Published online first, 2025.

\bibitem[MPT]{mpt26}
M.~Mourgoglou, C.~Puliatti, and X.~Tolsa.
\newblock Failure of the ${L}^p$ {N}eumann problem in bi-{L}ipschitz images of
  a ball.
\newblock Preprint. 2026. arXiv:2611.24519.

\bibitem[MT24]{mt24}
M.~Mourgoglou and X.~Tolsa.
\newblock The regularity problem for the {L}aplace equation in rough domains.
\newblock {\em Duke Math. J.}, 173(9):1731--1837, 2024.

\bibitem[Pog]{pog1}
B.~Poggi.
\newblock Failure to slide: a brief note on the interplay between the
  {K}enig-{P}ipher condition and the absolute continuity of elliptic measures.
\newblock Preprint. 2019. arXiv: 1912.10115.

\bibitem[She07]{shen07}
Z.~Shen.
\newblock A relationship between the {D}irichlet and regularity problems for
  elliptic equations.
\newblock {\em Math. Res. Lett.}, 14(2):205--213, 2007.

\bibitem[Ste70]{stein70}
E.~M. Stein.
\newblock Singular integrals and differentiability properties of functions,
  volume No. 30 of {\em Princeton Mathematical Series}.
\newblock Princeton University Press, Princeton, NJ, 1970.

\bibitem[TT24]{tt24}
O.~Tapiola and X.~Tolsa.
\newblock Connectivity conditions and boundary {Poincar\'e} inequalities.
\newblock {\em Analysis \& PDE}, 17(5):1831--1870, 2024.

\bibitem[Ulma]{ulmer3}
M.~Ulmer.
\newblock Equivalence between solvability of the dirichlet and regularity
  problem under an $l^1$ carleson condition on $\partial_t a$.
\newblock Preprint. 2026. arXiv:2509.10328.

\bibitem[Ulmb]{ulmer2}
M.~Ulmer.
\newblock The regularity problem with a weaker condition on only the
  transversal direction.
\newblock Preprint. 2025. arXiv:2503.20453.

\bibitem[Ulmc]{ulmer1}
M.~Ulmer.
\newblock Solvability of the {D}irichlet problem using a weaker {C}arleson
  condition in the upper half plane.
\newblock Preprint. 2025. arXiv:2503.19106.

\bibitem[Ulm26]{ulmer26}
M.~Ulmer.
\newblock Solvability of the {D}irichlet problem for a new class of elliptic
  operators.
\newblock {\em J. Differential Equations}, 485:Paper No. 114740, 74, 2026.

\end{thebibliography}

\end{document}